\documentclass{article}      
\usepackage{amssymb}
\usepackage{amsmath}
\usepackage{latexsym}
\usepackage[mathscr]{euscript}
\usepackage{graphicx}
\usepackage{amscd}
\usepackage{amsthm} 
\usepackage{fullpage} 
\usepackage{subfigure}
\usepackage{wrapfig}
\newtheorem{theorem}{Theorem}
\newtheorem{corollary}[theorem]{Corollary}
\newtheorem{lemma}[theorem]{Lemma}  
\newtheorem{proposition}[theorem]{Proposition}
\newtheorem{conjecture}[theorem]{Conjecture}
\newtheorem{definition}{Definition}

\newcommand{\bz}{\mathbb{Z}}
\newcommand{\br}{\mathbb{R}}

\newcommand{\p}{\partial}
\newcommand{\cc}{\mathcal{C}}

\newcommand{\cs}{\mathcal{S}}

\newcommand{\cp}{\mathcal{P}}

\newcommand{\cm}{\mathcal{M}}

\newcommand{\ca}{\mathcal{A}}

\newcommand{\hk}{\hookrightarrow}
\newcommand{\bg}{\bigskip}
\newcommand{\med}{\medskip}
\newcommand{\la}{\longrightarrow}
\newcommand{\bfl}{\begin{flushleft}}
\newcommand{\efl}{\end{flushleft}}

\newcommand{\eps}{\epsilon}
\newcommand{\G}{\Gamma}

\newcommand{\scrb}{\mathscr{B}}

\newcommand{\xr}{\xrightarrow}

 \newcommand{\ccg}{\cc_{\Gamma}}
 \newcommand{\ag}{Aut (\Gamma)}
 \newcommand{\ago}{Aut(\Gamma_1)}
 \newcommand{\agt}{Aut(\Gamma_2)}
 \newcommand{\ccgm}{\cm (\G, M)}
 \newcommand{\tccgm}{\tilde\cm (\G, M)}

\begin{document}  

\title{Morse theory, graphs, and string topology}
\author{Ralph L. Cohen   \thanks{The   author was partially supported by a grant from the NSF
}   \\
 Stanford University 
 }
   \date{\today}

 \maketitle

 \begin{abstract}
 In these lecture notes we discuss a body of work in which Morse theory is used to construct various homology and cohomology operations.  In the classical setting of algebraic topology    this is done by constructing a moduli space of graph flows, using homotopy theoretic methods to construct a virtual fundamental class, and evaluating
 cohomology classes on this fundamental class.  By using similar constructions based on ``fat" or ribbon graphs, we  describe how to construct
  string topology operations on the loop space of a manifold, using Morse theoretic techniques.   Finally,   we discuss how to relate these string topology operations to the counting of $J$- holomorphic curves in the cotangent bundle.  We end with speculations about the relationship between the absolute and relative Gromov-Witten 
 theory of the cotangent bundle,  and the open-closed  string topology of the underlying manifold.  
\end{abstract}

 \tableofcontents

 \section*{Introduction}  

Recently, an intersection theory has been developed for spaces of paths and loops in compact, oriented, manifolds.  This theory, which goes under the name of ``string topology",  was initiated in the seminal work of Chas and Sullivan \cite{chassullivan}, and has been expanded by several authors \cite{cohenjones}, \cite{cohengodin}, \cite{voronov}, \cite{chataur}, \cite{klein}, \cite{pohu}.   Various operad and field theoretic properties of this theory are now known \cite{voronov}, \cite{cohengodin}.  In particular, it was shown in \cite{cohengodin} how, given a surface $\Sigma$  viewed as a cobordism between $p$ circles and $q$-circles, there is an associated  operation in homology (as well as generalized homologies)
$$
\mu_\Sigma : H_*(LM)^{\otimes p} \to H_*(LM)^{\otimes q}.
$$
These operations share many formal properties with Gromov-Witten theory for
symplectic manifolds, and it is the  goal of an ongoing project of the author to understand
this relationship.

The possibilities for relating these theories became more compelling with the recent work
of Viterbo \cite{viterbo},  Salamon-Weber \cite{salamonweber},  and Abbodandolo-Schwarz \cite{abschwarz},  which proves that the Floer homology of the cotangent bundle of a closed, oriented manifold is isomorphic to the homology of the loop space,
$$
HF_*(T^*M) \cong H_*(LM).
$$
Here the Floer homology is defined  with respect to  a Hamiltonian $H_V$ on $T^*(M)$ defined in terms of a potential function on the manifold,  $V : \br / \bz \times M \to \br$.  

\med
Because Floer homology is a Morse homology based on the symplectic action functional, in order
to compare Floer theory and Gromov - Witten theory to string topology,  a first step  is to devise a Morse theoretic approach to string topology.  This will be the main topic of these lecture notes. 

The basic strategy is to expand and generalize a theory of graph flows of the author and Betz \cite{betzcohen} for constructing classical (co)homology invariants in algebraic topology.  Generalizing it in the appropriate way is the subject of joint work with P. Norbury, which we will report on in this paper.
The basic idea is to make ``toy model" of Gromov-Witten theory in which  the role of surfaces are replaced by finite graphs.  More specifically, we develop a theory in which we make the following
replacements from classical Gromov-Witten theory.   
\begin{enumerate}
\item A smooth surface $F$ is replaced by a finite, oriented graph $\Gamma$.
\item The role of the genus of $F$ is replaced by the first Betti number, $b = b_1 (\G)$.
\item The role of marked points in $F$ is replaced by the univalent vertices (or ``leaves") of $\G$.
\item The role of a complex structure $\Sigma$ on $F$ is replace by a metric on $\G$.
\item The notion of a $J$-holomorphic map to a symplectic manifold with compatible almost complex structure,  $\Sigma \to (N, \omega, J)$,  is replaced by the notion of a ``graph flow"  $\gamma : \G \to M$
which, when restricted to each edge is a gradient trajectory of a Morse function on $M$.  More specifically, each edge of the graph is labeled by a Morse function, $f_i : M \to \br$, and the restriction of $\gamma$ to the $ith$ edge  is a gradient trajectory of   $f_i$.  
\end{enumerate}

We   describe the moduli space of graph flows $\cm (\G, M)$, and then  define, using homotopy theory, virtual fundamental classes.  These are  \sl integral \rm homology classes in our case, and indeed we define these fundamental classes in any generalized homology theory that  supports an orientation of $M$.  In particular   we can therefore define  integral Gromov-Witten type invariants in this toy model situation.  These are   operations
$$
q_\G : H_*(BAut(\G)) \otimes H_*^{\Sigma_p} (M^p) \to H_*^{\Sigma_q}(M^q).
$$
Here $H_*^{\Sigma_n}(M^n)$ is the $\Sigma_n$-equivariant homology of the $n$-fold cartesian product $M^n$, where the symmetric group $\Sigma_n$ acts by permuting the coordinates.  $Aut(\G)$ is the finite group of combinatorial automorphisms of the graph $\Gamma$, and $BAut(G)$ is its classifying space.  The cohomology of this automorphism group is playing the role in this theory, of the cohomology of mapping class group (or equivalently the moduli space of Riemann surfaces) in Gromov-Witten theory. 

We review the results of \cite{betzcohen} and of \cite{cohennorbury} to show how one obtains classical cohomology operations such as cup product, intersection product, Steenrod squares, and Stiefel -Whitney classes  from these Morse-graph constructions.  We also discuss  certain field theoretic (or gluing) properties of these operations, as well as a kind of homotopy invariance. These properties are proved in \cite{cohennorbury}.  We then indicate
how these properties can be used to prove the classical relations, such as the Cartan and Adem relations, for these operations. 

Our next step is to move to the loop space.  Here we need to use ``fat" or ``ribbon graphs" instead of all finite graphs.  Spaces of fat graphs have been used by many authors to study moduli spaces of Riemann surfaces \cite{harer}, \cite{penner}, \cite{strebel}, \cite{kontsevich}.  We describe work of   Chas and Sullivan \cite{chassullivan}, Cohen and Jones \cite{cohenjones}, and Cohen and Godin \cite{cohengodin}, that show how to use fat graphs  to produce operations on the homology (or generalized homology)  of loop spaces,
$$
\mu_{g,n} : H_*(\cm_{g,n}) \otimes H_*(LM)^{\otimes p} \la H_*(LM)^{\otimes q}.
$$
Here $\cm_{g,n}$ is the moduli space of genus $g$-Riemann surfaces with $n = p+q$ parameterized boundary circles. These operations are constructed using classical techniques of algebraic and differential topology, such as transversal intersections of chains, and the Thom-Pontrjagin construction.  
 We also describe work of Sullivan \cite{sullivan2} and Ramirez \cite{ramirez} that gives an ``open-closed" version of these string topology operations.  This is motivated by open string theory in physics, where one studies configurations
of paths in a manifold  that take boundary values in ``D-branes".  

Next we adapt the Morse theoretic techniques described above to the loop space.
We describe work of the author \cite{cohen2} and Ramirez \cite{ramirez} that describe how to constuct the string topology operations using Morse theoretic techniques similar
to those used to construct classical (co)homology operations on finite dimensional manifolds.

Finally we review the work of  Salamon and Weber \cite{salamonweber} on the Floer homology of the cotangent bundle of a closed, oriented manifold, with its canonical symplectic structure.  They establish a relationship between moduli spaces of $J$-holomorphic cylinders in the cotangent bundle, and gradient flow lines of perturbed energy functionals on the loop space.  In both cases, they are using
pertubations by a potential function, $V : \br / \bz \times M \to \br$.  In the first case they use this to define a Hamiltonian on the cotangent bundle, $H_V : \br / \bz \times T^*M \to \br$, and in the second case
to define a Morse function $S_V : LM \to \br$.  We use this correspondence to relate
the space of gradient graph flows to the loop space, using a fat graph $\G$,  to the space of ``cylindrical $J$-holomorphic curves"  from a surface built from the graph $\G$  to the cotangent bundle, $T^*M$. In the case when the graph is a figure eight, we show how combining recent work of   Abbondandalo-Schwarz  \cite{abschwarz2} and Ramirez \cite{ramirez} implies that pair of pants product in the Floer homology of the cotangent bundle is equal to the Chas-Sullivan loop product on    $H_*(LM)$.

\med
The organization of the paper is the following.  In section 1 we recall results from \cite{betzcohen}
and describe results from \cite{cohennorbury} concerning the classical cohomology operations obtained using Morse theory, from  our ``toy model" of Gromov-Witten theory.  In section 2 we recall some of the theory of fat graphs, and describe the basic
constructions in string topology from \cite{chassullivan}, \cite{cohenjones}, \cite{cohengodin}, \cite{sullivan2},  and \cite{ramirez}.  In Section 3 we apply the methodology of studying graph flows to describe a Morse theoretic description of the (closed) string topology operations.  In section 4 we recall the results of Salamon-Weber \cite{salamonweber} and describe how they  might  be applied to get descriptions of string topology operations using ``cylindrical $J$-holomorphic maps" in the cotangent bundle.  We also   speculate about the further relationships between the Gromov-Witten theory (including the relative theory) of the cotangent bundle and the string topology of the manifold. 

\med
This paper is a survey of constructions and results of other papers, some of which have yet to appear.
I would like to thank my collaborators and students, J. Jones, V. Godin, P. Norbury, and A. Ramirez, as well as my colleague Y. Eliashberg, for many stimulating conversations concerning this and related material.  I would also like to thank O. Cornea , F. Lalonde, and the staff of the mathematics department at the University of Montreal for organizing such a lively and informative workshop.

\section{Graphs, Morse theory, and cohomology operations}

In this section we describe a generalization of the constructions and results of \cite{betzcohen}
that is joint work with P. Norbury \cite{cohennorbury}.  As described in the introduction, our goal
is to make a toy model of Gromov-Witten theory, where surfaces are replaced by graphs, and so on. 
The analogue of Teichmuller space in this theory is the space of graphs with metrics.  To describe
this we use ideas of Culler-Vogtman  \cite{cullervogtman}, and modifications of them due to Igusa \cite{igusa},  and to Godin \cite{godin}.  We define this space of structures in terms of a category of graphs.

\med
\begin{definition}\label{categ} Define $\cc_{b, p+q}$ to be the category of oriented graphs of first betti number $b$, with $p+q$ leaves.  More specifically, the objects of $\cc_{b, p+q}$ are finite graphs (one dimensional CW-complexes) $\G$,  with the following properties:
\begin{enumerate}
\item Each edge of the graph $\G$ has an orientation.
\item  $\G$ has  $p+q$ univalent vertices, or ``leaves".   $p$ of these are vertices of edges whose orientation points away from the vertex
(toward the body of the graph).  These are called  ``incoming"  leaves.  The remaining $q$ leaves are on edges whose orientation points
toward the vertex (away from the body of the graph).  These are called ``outgoing" leaves. 
\item $\G$ comes equipped with a  ``basepoint", which is a nonunivalent vertex.
\end{enumerate}
For set theoretic reasons we also assume that the objects in this category (the graphs) are subspaces  of a fixed infinite dimensional Euclidean space, $\br^\infty$. 

A morphism between objects $\phi : \G_1 \to \G_2$ is combinatorial map of graphs (cellular map) that satisfies:
\begin{enumerate}
\item $\phi$ preserves the orientations of each edge.
\item The inverse image of each vertex is a tree (i.e a contractible subgraph).
\item The inverse image of each open edge is an open edge.
\item $\phi$ preserves the basepoints.
\end{enumerate}
\end{definition}

Notice these conditions on morphisms is equivalent to saying that $\phi : \Gamma_1 \to \G_2$ is a combinatorial map of graphs that
is orientation preserving on edges and is a basepoint preserving homotopy equivalence on the geometric realizations of the graphs. 
\med
\begin{figure}[ht]
  \centering
  \includegraphics{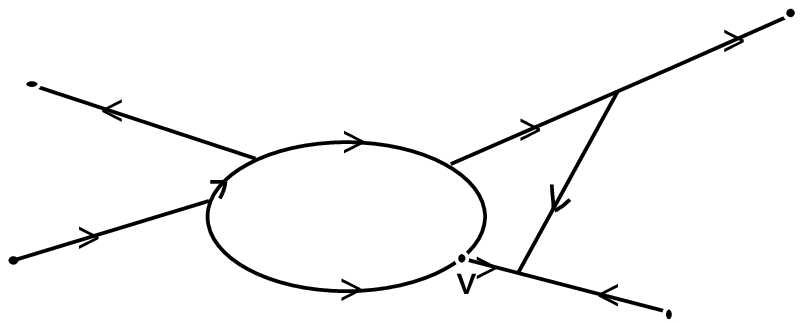}
  \caption{An object $\Gamma$ in $\cc_{2, 2+2}$}
  \label{fig:figone}
\end{figure}
\med
Given a graph $\Gamma \in \cc_{b, p+q}$, we define the automorphism group $Aut (\G)$ to be the group of invertible morphisms from 
$\Gamma$ to itself in this category.  Notice that $Aut (\G)$ is a finite group, as it is a subgroup of the group of permutations of the the edges. 

\med
We now fix a graph $\G$ (an object in $ \cc_{b, p+q}$), and we discuss the category of  ``graphs over $\G$", $\ccg$.  As we will see below, this category will be viewed
as the space of metrics on (subdivisions) of $\G$. 

\med
\begin{definition} Define $\ccg$ to be the category whose  objects are morphisms in $\cc_{b, p+q}$ with target $\G$:  $\phi : \G_0 \to \G$.
A morphism from $\phi_0 : \G_0 \to \G$ to $\phi_1 : \G_1 \to \G$ is a morphism $\psi : \G_0 \to \G_1$ in $ \cc_{b, p+q}$  with the property
that $ \phi_0  = \phi_1 \circ \psi  : \Gamma_0 \to \Gamma_1 \to \G$. 
\end{definition}

\med
We notice that the identity map $id : \G \to \G$ is a terminal object in $\ccg$.  That is, every object $\phi : \G_0 \to \G$ has a unique 
morphism to $id : \G \to \G$.  This implies that the geometric realization of the category, $|\ccg|$ is contractible.  ($|\ccg|$ is essentially the cone on the vertex represented by the terminal object.)   But notice that the category $\ccg$ has a free right action of the automorphism group, $Aut (\G)$,  given on the objects by composition:
\begin{align}
Objects \, (\ccg) \times Aut (\G) & \to Objects \, (\ccg) \notag \\
(\phi : \G_0 \to \G ) \cdot g &\to g \circ \phi : \G_0 \xrightarrow{\phi} \G \xrightarrow{g} \G
\end{align}

This induces an action on the geometric realization $\ccg$.  We therefore have:

\begin{proposition}\label{baut}
The orbit space is homotopy equivalent to the classifying space,
$$|\ccg|/ \ag  \simeq B\ag. $$
\end{proposition}

\med
As mentioned above, we can think of $|\ccg|$ as the space of ``metrics on subdivisions of $\G$".  The following idea of Igusa \cite{igusa}
associates to a point in $|\ccg|$ a metric on a graph over $\G$.
 
Recall that 
$$
|\ccg| = \bigcup_k \Delta^k \times
  \{\G_k \xr{\psi_k} \G_{k-1} 
  \xr{\psi_{k-1}}\G_{k-2} \to \cdots \xr{\psi_1} \Gamma_0 \xr{\phi}\G \} / \sim
 $$
where the identifications  come from the face and degeneracy operations. 

\begin{figure}[ht]
  \centering
  \includegraphics[height=12cm]{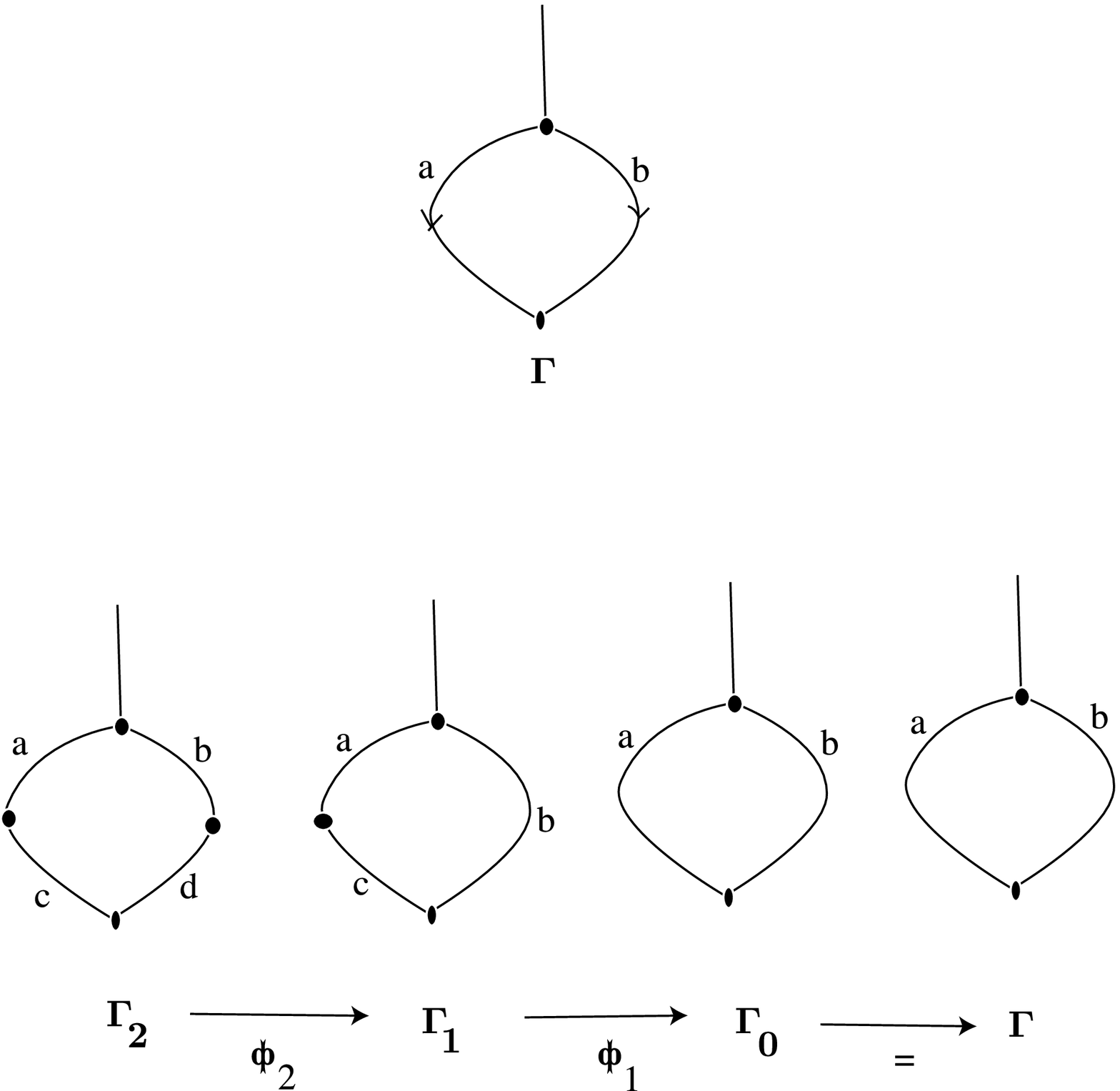}
  \caption{A $2$-simplex in $|\ccg|$.}
  \label{fig:figtwo}
\end{figure}

  Let $(\vec{t}, \vec{\psi})$ be a point in $|\ccg|$, where
$\vec{t} = (t_0, t_1, \cdots , t_k)$ is a vector of positive numbers whose sum equals one, and  $\vec{\psi}$ is a sequence of $k$-composable
morphisms in $\ccg$.   Recall that a morphism $\phi_i : \G_i \to \G_{i-1}$ can only collapse trees, or perhaps compose such a collapse with an automorphism.  So in a sense, given a composition of morphisms, 
$$
\vec{\psi} : \G_k \to \cdots \to \G_0 \to \G
$$ 
$\G_k$ is a (generalized)  subdivision of $\G$, and in particular $\G$ is obtained from $\G_k$ by collapsing various edges.  We use the coordinates $\vec{t}$ of the simplex $\Delta^k$ to define a metric on $\G_k$ as follows.  For each edge $E $ of $\G_k$, define $k+1$ numbers, $\lambda_0(E), \cdots , \lambda_k (E)$ given by
$$
\lambda_i(E) = \begin{cases} 0 \quad \text{if $E$ is collapsed by $\vec{\psi}$ in $\G_i$, and}, \\
1 \quad \text{if $E$ is not collapsed in $\G_i$} \end{cases}
$$

We then define the length of the edge $E$ to be
\begin{equation}\label{length}
\ell (E) = \sum_{i=0}^k t_i \lambda_i (E).
\end{equation}

\begin{figure}[ht]
  \centering
  \includegraphics[height=12cm]{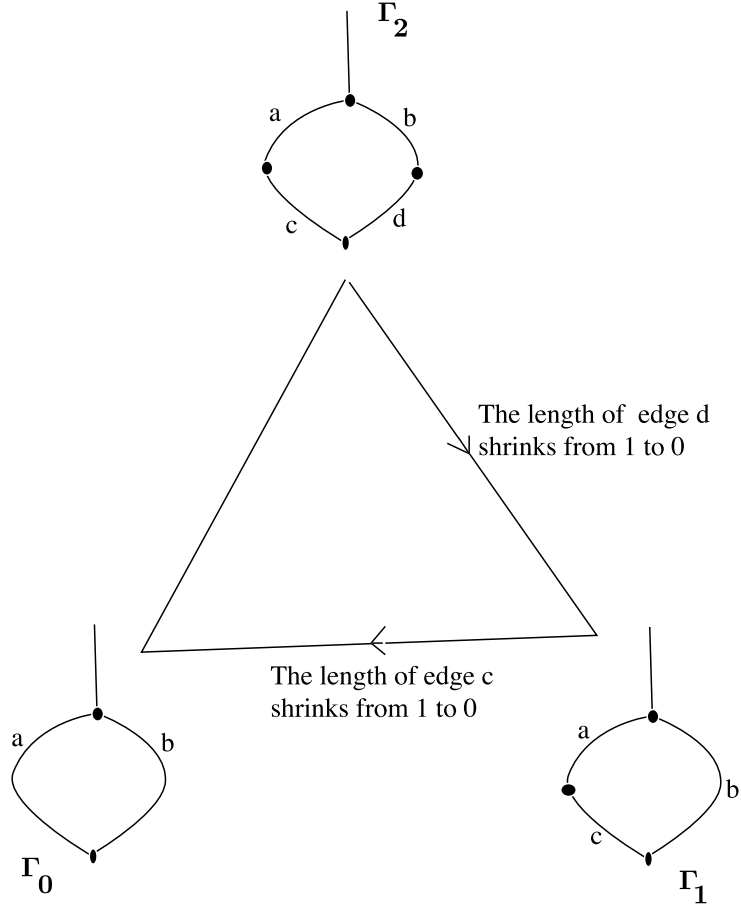}
  \caption{A $2$-simplex of metrics.}
  \label{fig:figthree}
\end{figure}

Notice also that  
the orientation
on the edges and the metric deterimine parameterizations (isometries) of standard intervals to the edges of the graph $\G_k$ over $\G$,
\begin{equation}\label{param}
\theta_E : [0, \ell (E)] \xr{\cong} E 
\end{equation}

Thus a point $ (\vec{t}, \vec{\psi}) \in|\ccg|$  determines a metric on a graph  $\G_k$ living over $\G$, as well as a parameterization of its edges. 

\med
 We now fix a target closed manifold $M$ of dimension $d$.  Our goal is describe a moduli space of maps from  the space of graphs  over $\G$ to $M$,     which, when restricted to each edge, satisfy a gradient flow equation of a function on $M$.  This moduli space of  maps (which we call ``graph flows") will be our analogue in this model, of the moduli space of $J$- holomorphic curves in a symplectic manifold.     For these purposes, we will need a bit more structure on an object,      $\phi_0 : \G_0 \to \G$.  Namely,  we want to label
the edges of $\G_0$ by distinct functions, $f_i : M \to \br$.  To  say this more precisely, we first introduce some notation.

Let $V$ be a  real vector space, and let $F(V, m)\subset V^m$ be the configuration space of $m$- distinct vectors in $V$.  It is an open, dense subset of the $m$-fold cartesian product, $V^m$.   We note that if $V$ is infinite dimensional, $F(V,m)$ is contractible.  Consider the functor 
$$
\mu : \ccg \to \, Spaces
$$
which assigns to a graph over $\G$,    $\phi_0 : \G_0 \to \G$, the configuration space $F(C^\infty (M), e(\G_0))$, where $C^\infty (M)$ is the vector space of smooth, real valued functions on $M$, and $e(\G_0)$ is the number of edges of $\G_0$.  Given a morphism
$\psi : \G_1 \to \G_0$, which collapses certain edges  and perhaps permutes others, there is an obvious induced map,
$$
\mu (\psi) : F(C^\infty (M),  e(\G_1)) \to F(C^\infty (M), e(\G_0)).
$$
This map projects off of the coordinates corresponding to edges collapsed by $\psi$, and permutes coordinates corresponding to  the permutation of edges induced by $\psi$.  

We can now do a homotopy theoretic construction, called the homotopy colimit (see for example \cite{bousfieldkan}).

\begin{definition}\label{structure} We define the space of metric structures and Morse labelings on $G$, $ \cs (\G)$,  to be the homology colimit,
$$\cs(\G) = hocolim \, (\mu : \ccg \to \, Spaces).$$
\end{definition}

The homotopy colimit construction is a simplicial space whose $k$ simplices consist of pairs,
$(\vec{f}, \vec{\psi})$, where $\vec{\psi} : \Gamma_k \to \Gamma_{k-1} \to \cdots \Gamma_0 \to \Gamma$ is a $k$-tuple of composable morphisms in $\ccg$, and $\vec{f} \in \mu (\Gamma_k)$.  That is, $\vec{f}$  is a labeling of the edges of $\G_k$ by distinct functions $f_E : M \to \br$, for $E$ an edge of $\G_k$.  We refer to $\vec{f}$ as a ``$M$ - Morse labeling" of the edges of $\G_k$. 

So we can think of a point $\sigma \in \cs  (\G)$ as consisting  of a metric on a graph over $\G$, together with an $M$-  Morse labeling of its edges.

We now make the following observation.

\begin{lemma}  The space of structures  $\cs (\G)$ is contractible with a free $\ag$ action.
\end{lemma}

\begin{proof}   The contractibility follows from  standard facts about the homotopy colimit construction,  considering the fact that both $|\ccg|$ and $F(C^\infty (M), m)$ are contractible.  The free action of $\ag$ on $|\ccg|$ extends to an action on $\cs (\G)$, since $\ag$ acts by permuting the edges of $\G$, and therefore permutes the labels accordingly.  \end{proof}

\med
We now define our moduli space of structures.

\med
\begin{definition}  The moduli space of metric structures and $M$-Morse labelings on $G$, $\cm (\G)$,  is defined to be the quotient,
$$
\cm (\G) = \cs (\G) / \ag.
$$
\end{definition}

\med 
We therefore have the following.

\begin{corollary} The moduli space is a classifying space of the automorphism group,
$$
\cm (\G) \simeq B\ag.
$$
\end{corollary}

\med
We now define the moduli space of ``graph flows" in $M$.   Let $\sigma \in \cs (\G)$.  As mentioned above, $\sigma$ can be thought
of as a graph $\G_k$ over $\G$, with a metric and parameterizations of its edges, and a labeling of its edges by functions on $M$.  
Let $\gamma : \G_k \to M$ be a continuous map which is smooth on the open edges.  If we restrict $\gamma$ to the edge $E$, and use the parameterization induced by the metric, this defines a smooth map 
$$
\gamma_E : [0, \ell (E)] \to M
$$
for each edge $E$ of $\G_k$.

\begin{definition}\label{graphflow} Define the structure space of ``graph flows" on $\G$, $\tccgm$, to be 
\begin{align}
\tccgm &= \{(\sigma, \gamma) \, : \,  \sigma \in \cs (\G), \, \gamma : \G_k \to M \, \text{ is such that  } \, \gamma_E : [0, \ell (E)] \to M  \notag \\ &\text{satisfies the gradient flow   equation}  \quad
 \frac{d\gamma_E}{dt} + \nabla f_E (\gamma (t)) = 0, \quad \text{for each edge $E$ of $\G_k$.} \} \notag
\end{align}
We define the moduli space of ``graph flows" on $\G$ $\ccgm$ to be the quotient,
$$
\ccgm = \tccgm/\ag.
$$
\end{definition}

\med
We discuss the topology of the moduli space of graph flows in \cite{cohennorbury}.  In any case the topology can be deduced from the following two results.

\med
\begin{theorem}\label{tree}  Suppose a graph $\G$ in $\cc_{b, p+q}$ is a tree.  Then there is a homeomorphism,
\begin{align}
\Psi : \ccgm &\xr{\cong} \cm (\G) \times M \simeq B\ag \times M \notag \\
(\sigma, \gamma)  &\la   \sigma \times \gamma (v)  \notag
\end{align}
where $v$ is the fixed vertex of the graph $\G_k$ over $\G$ determined by the structure $\sigma$. 
\end{theorem}
\begin{proof}  This follows from the existence and uniqueness theorem for solutions of ODE's on compact
manifolds.  The point is that the values of $\gamma$ on the edges emanating from $v$ are completely determined
by $\gamma (v) \in M$, since one has a unique flow line through that point for any of the functions labeling these edges.
The value of $\gamma$ on these edges determines the value of $\gamma$ on coincident edges (i.e edges that share a vertex)
for the same reason.  The theorem then follows.
\end{proof}

For general graphs $\G$, we analyze the topology of $\ccgm$ in the following way.
Let $\sigma \in \cs (\G)$.  A \sl tree  flow \rm of  $\Gamma$ with respect to the structure $\sigma$ is a collection $\gamma = \{\gamma_T\}$ where $\gamma_T : T \to M$ is a graph flow on a maximal subtree $T \subset \G_k$.  The collection ranges over all maximal subtrees $T \subset \G_k$, and is subject only to the condition that
the values at the basepoint are the same:
$$
\gamma_{T_1} (v) = \gamma_{T_2}(v)
$$
for any two maximal trees $T_1, T_2 \subset \G_k$.

We define
\begin{equation}
\cm_{tree}(\G, M) = \{(\sigma, \gamma) \, : \,  \sigma \in \cs (\G), \, \text{and} \,  
 \gamma = \{\gamma_T\} \, \text{is a tree flow of } \, \Gamma \, \text{with respect to} \, \sigma \}/ \ag
\end{equation}

We now have the following generalization of   theorem \ref{tree}.

\begin{theorem}\label{embed} 1.   For any graph  $\G \in \cc_{b, p+q}$ there is a homeomorphism
\begin{align}
\Psi : \cm_{tree}(\G, M) &\xr{\cong} \cm (\G) \times M \simeq B\ag \times M \notag \\
(\sigma, \gamma)  &\la   \sigma \times \gamma (v)  \notag
\end{align}
where $v$ is the fixed vertex of the graph $\G_k$ over $\G$ determined by the structure $\sigma$. 

2.  Let  $\rho : \ccgm \to \cm_{tree}(\G, M)$ be the map that sends a graph flow $\gamma$ to the tree flow obtained by restricting $\gamma$ to each maximal tree. Then
$\rho$ is a codimension $b_1(\Gamma)\cdot d$ embedding.
\end{theorem}

\med
\sl Idea of proof.  \rm The first part of this theorem follows from the existence and uniqueness theorem, just like theorem \ref{tree}.  For the second part, we will show in \cite{cohennorbury} that there exists a Serre fibration over the cartesian product,
$$
ev : \cm_{tree}(\G, M) \to M^{2b}
$$
where $b = b_1(\G)$, so that the restriction to the diagonal embedding
$$
\Delta^b : M^b \hk (M^2)^{b}
$$
is equal to $\ccgm$.  That is, there is a pullback square of fibrations,
\begin{equation}\label{evaluate}
\begin{CD}
\ccgm @>\rho >\hk >   \cm_{tree}(\G, M) \\
@V ev VV   @VV ev V \\
M^b @>\hk >\Delta^b > (M^2)^b
\end{CD}
\end{equation}

To give an idea of the map $ev : \cm_{tree}(\G, M) \to M^{2b}$, assume for simplicity that $b_1(\G) = 1$.  Let $\sigma \in \cs (\G)$ be a structure.  Then any maximal tree in $\G_k$  is obtained  by removing a single edge from $\G_k$.  
Say $T_E  $ is a maximal tree obtained by removing the edge $E$ of $\G_k$.
 Let  $\gamma_{T_E} : T_E \to M$ be a graph flow on $T_E$ with respect to $\sigma$. $E$ is an oriented edge, so we can identify its vertices as a source vertex $v_0$ and a target
 vertex $v_1$.   Let $ev_1(\gamma_{T_E}) \in M$ be the evaluation $\gamma_{T_E}(v_1)$.  Now 
 $\gamma_{T_E}$ is not defined on $E$, but we can extend $\gamma_{T_E}$ to the edge
 $E$ by considering the unique gradient flow line $\alpha_E : E \to M$ of the function
 $f_E :M \to \br$   labeling $E$ given by   the structure $\sigma$, that has the property
 that $\alpha_E(v_0) = \gamma_E(v_0)$.  We then define $ev_2(\gamma_E) \in M$
 to be the evaluation $\alpha_E(v_1)$.  Taking the pair $ev_1 \times ev_2$ defines
  an element
  $$ev (\gamma_{T_E}) = (\gamma_{T_E}(v_1), \alpha_E(v_1)) \in M^2.
  $$
  Clearly the graph flow $\gamma_{T_E}$ on the tree $T_E$ is the restriction of a graph flow on all of $\G_k$ if and only if $\gamma_E(v_1) = \alpha_E(v_1)$.  This is the basic ingredient in the construction of the fibration $ev : \cm_{tree}(\G, M) \to (M^2)^b$, and verifying the pullback property (\ref{evaluate}).   Details will appear in \cite{cohennorbury}.
  
  \med
  This result will be used to define virtual fundamental classes for these moduli spaces.  This will depend on the following
  homotopy theoretic construction described in \cite{hatcherquinn}, \cite{cohenjones}, and in full generality in 
  \cite{klein2}.
  
  \med
  \begin{proposition}\label{thom}
  Let $\iota : P \hk N$ be a smooth embedding of closed manifolds with normal bundle $\nu (\iota)$.  Let $p : E \to N$ be a Serre fibration, $\iota^*(E) \to P$
  the restriction to $P$.  Let $(\iota^*(E))^{\nu (\iota)}$ be the Thom space of the pullback of the normal bundle  via the projection map
  $p : \iota^*(E) \to P$.  Then there is a   ``Thom collapse map"
  $$
  \tau : E \to ( \iota^*(E))^{\nu (\iota)}
  $$
  compatible with the usual Thom collapse of the embedding, $N \to P^{\nu(\iota)}.$
  \end{proposition}
  
  \med
  Notice that there is no smoothness requirement about the fibration $E \to N$.  In particular this proposition implies that if $h_*$ is any (generalized)  homology theory for which the normal bundle
  $\nu (\iota)$ is orientable then one can define an ``umkehr map",
  $$
  \iota_{!} : h_q(E) \to h_q((\iota^*(E))^{\nu (\iota)}) \xr{\cong} h_{q-k}(\iota^*(E)),
  $$
  where $k$ is the codimension of the embedding $\iota$, and the second map in this composition is the Thom isomorphism. (The existence of the Thom isomorphism is exactly what is meant by saying the normal bundle is orientable with respect to the homology theory $h_*$.)   Equivalently, (or dually), one gets a ``pushforward map"  in cohomology,
  $$
  \iota^! : h^r(\iota^*(E)) \to h^{r+k}(E).
  $$
  
  So by theorem \ref{embed}, and in particular the structure displayed in diagram
  \ref{evaluate},  we can therefore construct a Thom collapse
  $$
  \tau : \cm (\G) \times M \to (\ccgm)^\nu
  $$
  where $\nu$ is the pullback along the evaluation map $ev : \ccgm \to M^b$
  of the normal bundle to the diagonal embedding $\Delta^b : M^b \hk (M^2)^{b},$
  which is isomorphic to the exterior product of the tangent bundle,
  $$
  \nu= (TM)^b \to M^b.
  $$
  Thus if $h_*$ is any homology theory that supports an orientation of $M$ (i.e the tangent bundle $TM)$,  there is an  umkehr map,
  $$
  \rho_! : h_*(B\ag \times M) \cong h_*(\cm (\G) \times M) \to h_{*-b\cdot dim (M)}(\ccgm).
  $$

   Now notice that since $\ag$ is a finite group, its classifying space $B\ag$ is infinite dimensional. Thus
  by theorem \ref{embed} the moduli space $\ccgm$ is infinite dimensional.  So in order to construct a fundamental
  class, we need to ``cut down" the    moduli  space $\cm (\G) \simeq B\ag$  by a homology class.  When we do that we can make the following definition.
  
  \begin{definition}\label{virtual} Let $h_*$ be a generalized homology theory that supports an orientation of $M$. 
   Given a homology class $\alpha$ in  $h_k(B\ag) = h_k(\cm (\G))$, we can  define a  \sl virtual fundamental class \rm
   $$
   [\cm^\alpha(\G,M)] \in h_{k+ \chi (\G) \cdot dim(M)}(\ccgm)
   $$
   by the formula 
   $$
    [\cm^\alpha(\G,M)] = \rho_! (\alpha \times [M])
    $$
    where $[M] \in h_{dim (M)}(M)$ is the fundamental class.
    \end{definition}
    
    \med
 
    Although  using generalized homology theories such as $K$-theory and bordism theory are very useful, and will be pursued in future work, for the rest of this paper we restrict to ordinary homology with integral coefficients, or perhaps with $\bz /2$ coefficients if the manifolds in question are not orientable. 
    
    \med

    \bf Geometric explanation.  \rm  The idea for defining this virtual fundamental class is the following. Suppose
    that the homology class $\alpha \in H_k(\cm (\G))$ is represented by a manifold $N_\alpha \subset \cm (\G)$.  
    Let $\tilde N_\alpha \subset \cm (\G))$ be the induced $\ag$- covering.   We can then define a subspace
    $$
    \tilde \cm^{N_\alpha}(\G, M) \subset \tilde \cm (\G, M)
    $$
    to consist of those pairs $\{(\sigma, \gamma) : \, \sigma \in \tilde N\}$.  We then define the ``cut down" moduli space
   $$\cm^{N_\alpha}(\G, M) =  \tilde \cm^{N_\alpha}(\G, M) / \ag.$$  Notice that the embedding $\rho$ restricts to give a codimension $b\cdot d$ embedding,
   $$
   \rho_\alpha : \cm^{N_\alpha}(\G, M) \hk N_\alpha \times M.
   $$
   Now \sl if \rm $N_\alpha $ could be chosen in such a way that $ \cm^{N_\alpha}(\G, M)$ is a \sl smooth, \rm
oriented  submanifold, which is compact, or could be compactified, then Poincare duality would imply that its fundamental class would be $\rho_!([N_\alpha] \times [M])$.  This class    would be independent of the particular manifold $N_\alpha$  we chose to represent $\alpha$, so long as the smoothness and compactness properties hold. Verifying that such manifolds \sl can \rm be chosen to satisfy these properties, would be a technical analytical problem.     By defining our virtual fundamental classes using the algebraic topological methods    represented in theorems (\ref{embed}) and (\ref{tree}), we avoid the analytical issues of smoothness and compactness in our construction.

 \med

  Given a graph flow $(\sigma, \gamma) \in \tccgm$, we can evaluate $\gamma$ on the $p+q$ univalent vertices (leaves) to define a map
\begin{equation}\label{eval}
  ev : \tccgm \la (M^p) \times (M^q).
\end{equation}
Notice furthermore  that if we restrict an automorphism in $\ag$  to the incoming and outgoing leaves,
there is a group homomorphsim
\begin{equation}\label{restrict}
r = r_{in} \times r_{out} : \ag \to \Sigma_p \times \Sigma_q.
\end{equation}
 The  evaluation map is clearly equivariant with respect to this homomorphism. This will allow
 us to pass to the (homotopy) orbit space to define an evaluation map
 \begin{equation}\label{eqeval}
 ev : \ccgm = \tccgm/\ag \la \left(E\Sigma_p \times_{\Sigma_p} M^p\right) \times \left(E\Sigma_q \times_{\Sigma_q} M^q\right)
 \end{equation}
 Here $E\Sigma_k$ is  a contractible space with a free $\Sigma_k$ action.  This type of passage to homotopy orbit spaces will be discussed in detail in \cite{cohennorbury}.
 
 By pulling back equivariant cohomology classes of $(M^p)$ and $(M^q)$ to $H^*(\ccgm)$ and then evaluating on the virtual fundamental classes defined above, then  for each $\alpha \in H_*(B\ag)$  we obtain an 
 operation 
\begin{align}
 q_\Gamma(\alpha) :  H^*_{\Sigma_p}(M^p) \otimes  H^*_{\Sigma_q}(M^q) &\la k \notag \\
 x \otimes y &\la \langle ev^*(x \times y),  [\cm^\alpha(\G,M)] \rangle,  
\end{align}
 where $k$ is any ground field taken as the coefficients of the cohomology groups. 
 
 Using the ``umkehr" maps described above, in \cite{cohennorbury}  we actually study the corresponding adjoint operations in equivariant homology,
\begin{equation}\label{operation}
 q_\G : H_*(B\ag) \otimes H_*^{\Sigma_p}(M^p) \la H_*^{\Sigma_q}(M^q).
\end{equation} 

In order to describe the properties of these operations, consider the following constructions.
\begin{enumerate}
\item Suppose $\G$ and $\G'$ are two graphs (objects) in the category $ \cc_{b, p+q}$, and 
$\phi : \G \to \G'$
a morphism.  By definition \ref{categ} this morphism induces a homomorphism between their automorphism groups and the associated classifying spaces,
$$
 \phi_* : Aut \G \la Aut \G' \quad B\phi : B Aut \G \la B Aut \G'.
$$
\item Let $\G_1$ be a graph  in $ \cc_{b, p+q}$, and $\G_2$ a graph in $\cc_{b', q+r}$.
Then let $\G_1 \# \G_2$ be the graph in $\cc_{b+b'+q-1, \, p+r}$ obtained by gluing the $q$ outgoing leaves of $\G_1$ to the $q$ incoming leaves of $\G_2$.
\begin{figure}[ht]
  \centering
  \includegraphics[height=10cm]{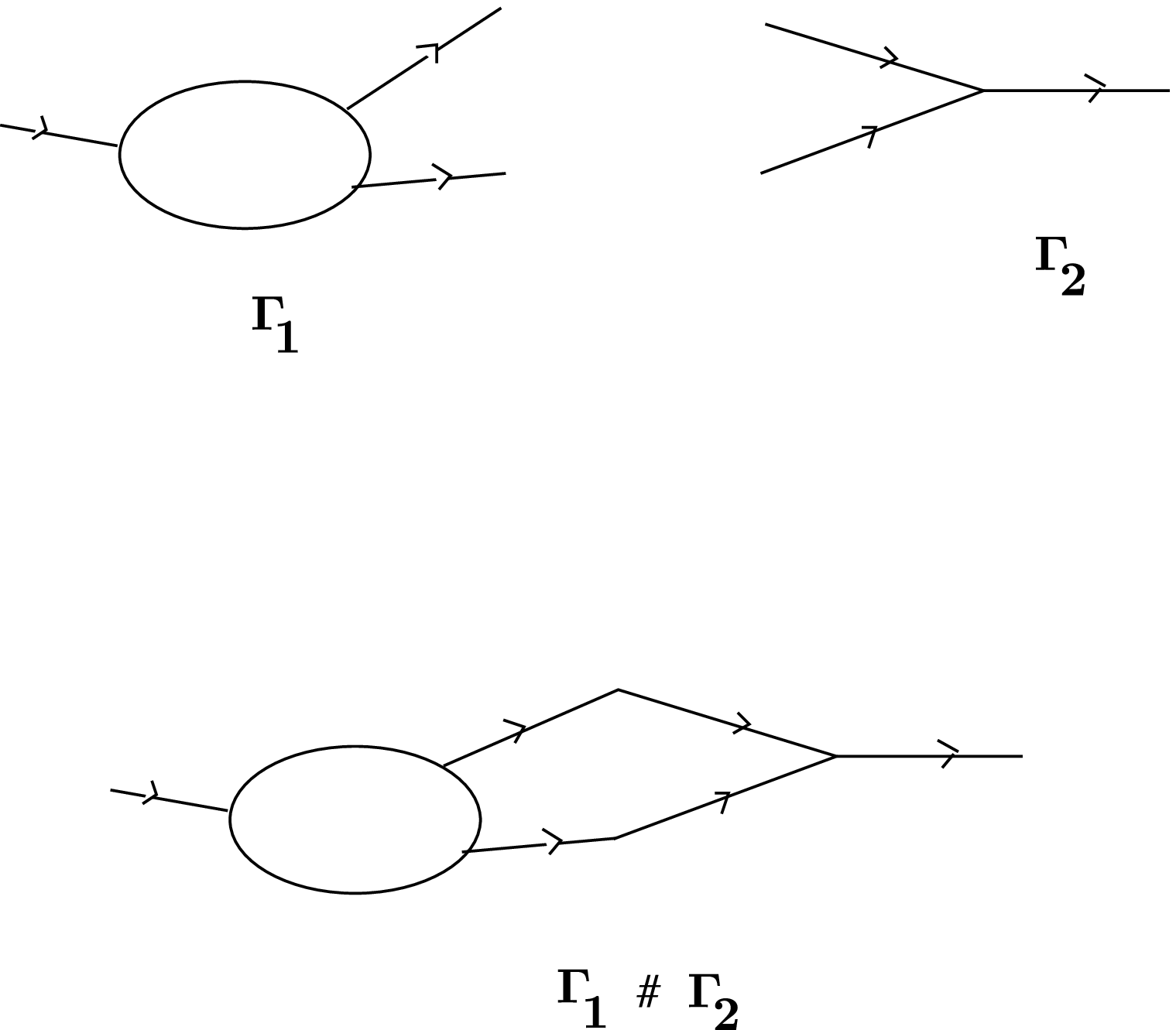}
  \caption{ $\G_1 \# \G_2$}
  \label{fig:figfour}
\end{figure}
\item   Let $\G_1$ and $\G_2$ be as above.  Consider the homomorphisms
$$\rho_{out} : \ago \to \Sigma_q \quad \rho_{in} : \agt \to \Sigma_q$$
  defined by the induced permutations of the outgoing and incoming leaves, respectively.   Let $\ago \times_{\Sigma_q} \agt$ be the fiber product of these homomorphisms.  That is, 
$$
\ago \times_{\Sigma_q} \agt \subset \ago \times \agt
$$
is the subgroup consisting of those $(g_1, g_2)$ with $\rho_{out}(g_1) = \rho_{in}(g_1)$.
Let 
$$
p_1 : \ago \times_{\Sigma_q} \agt \to \ago \quad \text{and} \quad p_2 : \ago \times_{\Sigma_q} \agt \to \agt
$$
be the projection maps.  There is also an obvious inclusion as a subgroup of the automorphism group of the glued graph, 
$$
\iota : \ago \times_{\Sigma_q} \agt  \hk Aut(\G_1 \# \G_2)
$$
which realizes  $ \ago \times_{\Sigma_q} \agt$ as the subgroup of $  Aut(\G_1 \# \G_2)$ consisting of automorphisms that preserve the subgraphs, $\G_1$ and $\G_2$.
\end{enumerate}

In \cite{cohennorbury} we prove the following theorem.

\med
\begin{theorem}\label{field}
The operations 
$$
q_\G : H_*(B\ag) \otimes H_*^{\Sigma_p}(M^p)  \to H_*^{\Sigma_q}(M^q)
$$
satisfy the following properties.
\begin{enumerate}
\item (Homotopy invariance)

 Let $\phi : \G \to \G'$ be a morphism in $ \cc_{b, p+q}.$ Then  the following diagram commutes:
 $$
 \begin{CD}
 H_*(B\ag) \otimes H_*^{\Sigma_p}(M^p) @>\phi_* \otimes 1>> H_*(B Aut(\G')) \otimes H_*^{\Sigma_p}(M^p) \\
 @Vq_{\G} VV     @VVq_{\G'} V \\
  H_*^{\Sigma_q}(M^q)  @>>= >   H_*^{\Sigma_q}(M^q).
  \end{CD}
  $$
  \item (Gluing formula)
  
  Let $\G_1 \in  \cc_{b, p+q}$, and $\G_2\in \cc_{b', q+r}$.  Let $\alpha \in H_*(B ( \ago \times_{\Sigma_q} \agt) ),$ and $x \in H_*^{\Sigma_p}(M^p)$.   Then
 $$ 
 q_{\G_1 \# \G_2}(\iota_*(\alpha) \otimes x) = q_{\G_2}\left(\left((p_2)_*(\alpha) \otimes  q_{\G_1}((p_1)_*(\alpha) \otimes x) \right)\right) \in H_*^{\Sigma_q}(M^q).
  $$
\end{enumerate}
\end{theorem}

\med
\bfl \bf Examples. \rm
\begin{enumerate}
\item   Consider the graph $\G = $
\begin{figure}[ht]
  \centering
  \includegraphics[height=4cm]{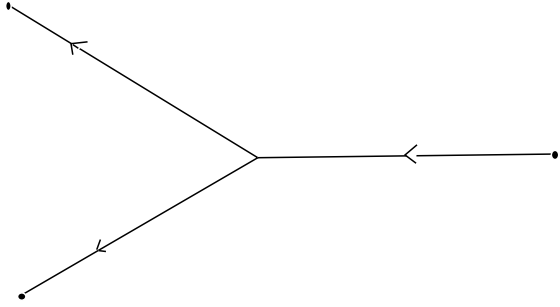}
\caption{The ``Y-graph"}
   \label{fig:figfive}
\end{figure}

$\G$ has one incoming and two outgoing leaves.  The automorphism group $\ag = \bz / 2$.  So the operation is a map
$$
q_{\G} : H_*(B\bz /2) \otimes H_*(M) \to H_*^{\Sigma_2}(M \times M).
$$
It is easy to see that this map is induced by Steenrod's ``equivariant diagonal"
$$
B \bz/2 \times M \to E\bz/2 \times_{\bz/2} (M\times M )
$$
given by applying the construction $E\bz/2 \times_{\bz/2} ( -)$ to the diagonal map
$\Delta : M \to M \times M$.   Now work with  $\bz/2$ coefficients. The Steenrod squares are defined if we take the dual
map in cohomology.  Namely, if $\alpha \in H^k(M)$, then  
$$
(q_{\G})^*(\alpha \otimes \alpha) = \sum_{i=0}^k a^k \otimes Sq^{k-i}(\alpha).
$$
Here $a \in H^1(B\bz/2 ; \bz/2) \cong \bz/2$ is the generator. 
It will be shown in \cite{cohennorbury} that the Cartan and Adem relations for the Steenrod squares follow from the homotopy invariance and gluing formulas in theorem \ref{field}. 

\item  Now consider the graph $\G $ given in figure 6 below. 
\begin{figure}[ht]
  \centering
  \includegraphics[scale=0.7]{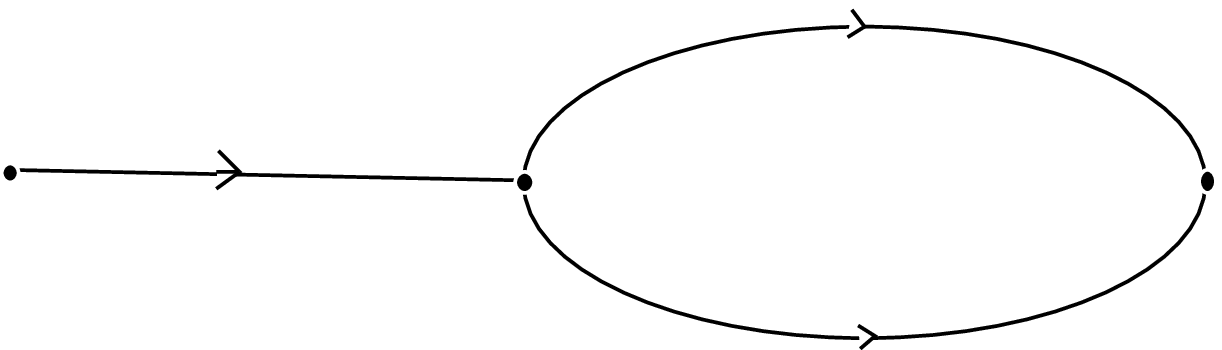}
\caption{}
   \label{fig:figsix}
\end{figure}

In this case there is only one incoming leaf, and the automorphism group is also $\bz/2$.  So the invariant (with $\bz/2$-coefficients)  is a map 
$$
q_{\G} : H_*(B\bz /2) \otimes H_*(M) \to \bz/2,
$$
or equivalently, $q_{\G} \in  H^*(B\bz /2) \otimes H^*(M).$
In \cite{cohennorbury} we will show   that
$$
q_{\G} = \sum_{i=0}^d a^i \otimes w_{d-i}(M)
$$
where $d = dim (M)$, and $w_j(M) \in H^j(M; \bz/2)$ is the $j^{th}$ Stiefel-Whitney class of the tangent bundle. 
\end{enumerate}

\efl 

\section{String topology}
In this section we review the basic constructions of   ``string topology".  This is an intersection theory for 
loop spaces and path spaces in manifolds, invented by Chas and Sullivan \cite{chassullivan}.  Their theory was constructed by studying chains in a manifold and in the loop space.  In our discussion we shall
follow the approach used in \cite{cohenjones} and \cite{cohengodin} based on the ``Thom collapse map".
This will allow us to identify a Morse theoretic approach to string topology based on graph operations, similar to those described in the last section.  This, in turn, will allow us to relate string topology to the counting of holomorphic curves in the cotangent bundle of the loop space.

To motivate the constructions of string topology, we first recall how classical intersection
is done in the differential topology of compact manifolds.

\med
Let $e : P^p \hk N^n$ be an embedding of smooth, closed, oriented manifolds. $p$ and $n$ are the dimensions of $P$ and $N$ respectively, and  let  $k$ be  the   codimension,   $k = n-p$.
The idea in intersection theory is to take a $q$-cycle  in $N$, which is transverse to $P$ in an appropriate sense, and take the intersection to produce a  $q-k$-cycle in $P$.  Homologically, one can make this rigorous by using Poincare duality, to define the intersection map,
$$
e_! : H_q(N) \to H_{q-k}(P)
$$
by the  composition
$$ 
e_! : H_q(N) \cong H^{n-q}(N) \xr{e^*} H^{n-q}(P) \cong H_{q-k}(P)
$$
where the first and last isomorphisms are given by Poincare duality.  

Intersection theory can also be realized by the ``Thom collapse" map.  Namely, extend the
embedding $e$ to a tubular neighborhood, $P \subset \eta_e \subset N$, and consider the projection map,
 $
\tau_e : N \to N/(N-\eta_e).
$
By the tubular neighborhood theorem, this space is homeomorphic to the Thom space of the normal bundle, $N/(N-\eta_e) \cong P^{\eta_e}$.  So the Thom collapse map can be viewed as a map,
$$
\tau_e : N \to P^{\eta_e}.
$$
Then the homology intersection map $e_!$ is equal to the composition,
\begin{equation}\label{umkehr}
e_! : H_q(N) \xr{(\tau_e)_*} H_q(P^{\eta_e}) \cong H_{q-k}(P)
\end{equation}
where the last isomorphism is given by the Thom isomorphism theorem.  In fact this description of the intersection (or ``umkehr") map $e_!$ shows that it can be defined in \sl any \rm generalized homology theory, for which there exists a Thom isomorphism for the normal bundle.  This is an orientation condition. 
In these notes we will restrict our attention to ordinary homology, but intersection theories in such (co)homology theories as $K$-theory and cobordism theory are very important as well. 

We remark that the ``intersection product" is the example induced by  the diagonal embedding,
$\Delta : M \to M\times M$.  This induces a product, 
\begin{equation}\label{intersect}
\Delta_! : H_p(M)\otimes H_q(M)) \to H_{p+q-d}(M),
\end{equation} where $d$ is the dimension of $M$.

\med 
The Chas-Sullivan ``loop product" in the homology of the free loop space of a closed oriented $n$-manifold,
\begin{equation}\label{loop}
\mu : H_p(LM) \otimes H_q(LM) \to H_{p+q-d}(LM)
\end{equation}
is defined as follows.  

Let $Map (8, M)$ be the mapping space from the figure 8  (i.e the wedge of two circles) to the manifold $M$.  The maps are required to be piecewise smooth (see \cite{cohenjones}).  
Notice that  $Map (8, M)$  is the  subspace of $LM \times LM$ consisting of those pairs of loops that agree at the  basepoint $1 \in S^1$.  In other words, there is a pullback square 
$$
\begin{CD}
 Map (8, M)   @>e >>  LM \times LM \\
 @V ev VV  @VV ev \times ev V \\
 M @>>\Delta > M \times M 
 \end{CD}
 $$
 where $ev : LM \to M$ is the fibration given by evaluating a loop at $1 \in S^1$.  The map $ev : Map (8, M)\to M$ evaluates the map at the crossing point of the figure 8.   Since $ev \times ev$ is a fibration,
$e : Map(8,M) \hk LM \times LM$ can be viewed as a codimension $d$ embedding, with normal bundle
$ev^*(\eta_\Delta) \cong ev^*(TM)$.   As was done in \cite{cohenjones} this diagram allows us to define a Thom-collapse map
$$
\tau_e : LM \times LM \to Map(8,M)^{ev^*(TM)},
$$
and therefore an intersection map,
$$
e_! : H_*(LM \times LM) \to H_{*-d}(Map(8,M)).
$$
Now by going around the outside of the figure 8, there is a map
$\gamma :  Map(8,M) \to LM$, and the Chas-Sullivan pairing is the composition,
$$
\mu_*  : H_p(LM) \otimes H_q(LM) \xr{e_!} H_{p+q-d}(Map(8,M)) \xr{\gamma_*} H_{p+q-d}(LM)
$$
which defines an associative, commutative algebra structure.  

\begin{figure}[ht]
  \centering
  \includegraphics[scale=0.7]{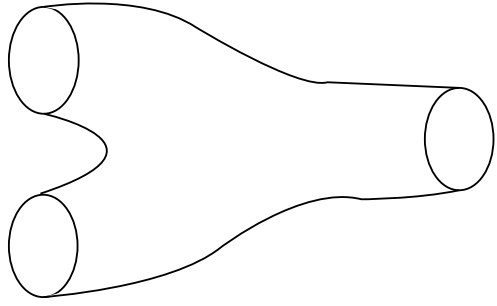}
\caption{The ``pair of pants" surface $P$}
   \label{fig:figseven}
\end{figure}

\med
One can think of this structure in the following way.  Consider the ``pair of pants" surface $P$, viewed as a cobordism from two circles to one circle (see figure \ref{fig:figseven}). 

Consider the smooth mapping space, $Map(P, M)$.  Then there are restriction
maps to the incoming and outgoing boundary circles,
$$
\rho_{in} : Map(P, M) \to LM \times LM, \quad \rho_{out} : Map(P,M) \to LM.
$$
Notice that the figure 8 is  homotopy equivalent to the surface $P$, with respect to which 
the restriction map $\rho_{in} : Map(P,M) \to LM \times LM$ is homotopic to the embedding
$e : Map (8,M) \to LM \times LM$.  Also restriction to the outgoing boundary, $\rho_{out} : Map(P,M) \to LM$
is homotopic to $\gamma : Map (8,M) \to LM$.  So the Chas-Sullivan product can be thought of as a composition,
$$
\mu_* :  H_p(LM) \otimes H_q(LM) \xr{(\rho_{in})_!} H_{p+q-d}(Map(P,M)) \xr{(\rho_{out})_*} H_{p+q-d}(LM).
$$
The role of the figure 8 can therefore be viewed as just a technical one, that allows us to define the umkehr map
$e_! = (\rho_{in})_!$.

More generally, consider a surface $F$, viewed as a cobordism from $p$-circles to $q$-circles.  See figure \ref{fig:figeight} below. 

\begin{figure}[ht]
  \centering
  \includegraphics[height=6cm]{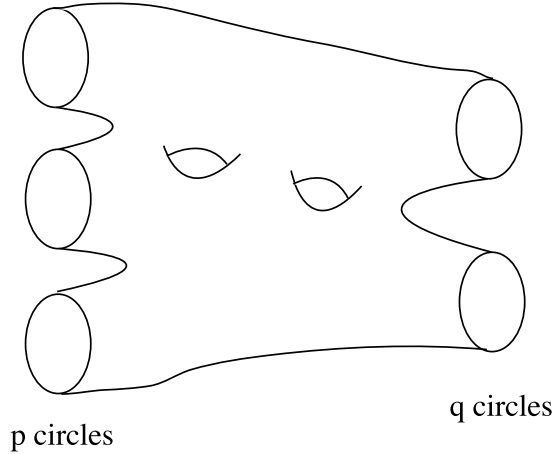}
\caption{The  surface $F$}
   \label{fig:figeight}
\end{figure}

We can consider the mapping space, $Map (F, M)$, and the resulting restriction maps,
\begin{equation}\label{restrict1}
(LM)^q \xleftarrow{\rho_{out}} Map (F, M) \xr{\rho_{in}} (LM)^p.
\end{equation}

The idea in \cite{cohengodin} is to construct an  ``umkehr map"
$$
(\rho_{in})_! : h_*((LM)^p) \to h_{* +\chi (F)\cdot d}(Map (F, M))
$$ where $\chi (F)$ is the Euler characteristic of the surface $F$ and $d = dim (M)$. Like above, $h_*$ is any generalized homology theory that supports an orientation of $M$.  This then allows
the definition of a string topology operation
$$
\mu_F : h_*((LM)^p) \xr{(\rho_{in})_!} h_{* +\chi (F)\cdot d}(Map (F, M)) \xr{(\rho_{in})_*} h_{*+\chi (F)\cdot d}((LM)^q).
$$

It was proved in \cite{cohengodin} that the operations $\mu_F$ respect  gluing of surfaces.
That is, if $F_{g,p+q}$ is a cobordism between $p$ circles and $q$-circles of genus $g$, and $F_{h,q+r}$ is a cobordism between $q$-circles and $r$ circles of genus $h$, and
$F_{g+h+q-1, p+r}$ is the glued cobordism between $p$ circles and $r$ circles, then  there is a    relation
$$
\mu_{F_{g+h+q-1, p+r}} = \mu_{F_{h,q+r}} \circ  \mu_{F_{g,p+q}} .
$$
This can be described in field theoretic language.  Notice in this case the number of outgoing
boundary components must be positive.  In \cite{cohengodin} this was condition was referred to as a ``positive boundary" condition, and the following was proved.

\med
\begin{theorem}\label{stringop}  If $h_*$ is a homology theory with respect to which the closed $d$-manifold $M$ is oriented, then the string topology operations which assigns to a circle $S^1$ the homology $h_*(LM)$,  and to a cobordism $F$, the operation $\mu_F : h_*((LM)^p) \to h_*((LM)^q)$ is a two dimensional, positive boundary, topological quantum field theory.
\end{theorem} 

\med
Clearly the difficult part in defining the string topology operations $\mu_F$ is the definition of the umkehr map $(\rho_{in})_!$. To do this, Cohen and Godin used the Chas-Sullivan idea of representing  the pair of pants  surface $P$ by a figure 8, and realized the surface $F$ by a ``fat graph" (or ribbon graph).   Fat graphs have been used to represent surfaces  for many years, and to great success.  See for example the following important works: \cite{harer}, \cite{strebel}, \cite{penner}, \cite{kontsevich}.  

We recall the definition.
\begin{definition}   A fat graph is a finite graph with the following properties:
\begin{enumerate} 
\item Each vertex is at least trivalent 
\item Each vertex comes equipped with a cyclic order of the half edges emanating from it.
\end{enumerate}
\end{definition}

We observe that the cyclic order of the half edges is quite important in this structure.
It allows for the graph to be ``thickened" to a surface with boundary.  This thickening, which will be defined carefully below, can be thought of as assigning a ``width" to the ink used in drawing a fat graph.  Thus one is actually drawing a two dimensional space, and it is not hard to see that it is homeomorphic to a smooth surface.  Consider the following two examples (figure 9) of fat graphs which consist of the same underlying graph, but have different cyclic orderings at the top vertex.

\begin{figure}[ht]
  \centering
  \includegraphics[height=8cm]{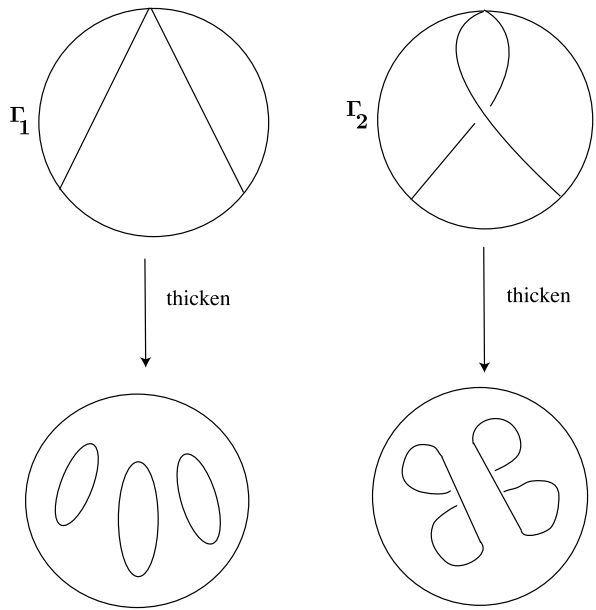}
 \caption{Thickenings of two fat graphs }
   \label{fig:fignine}
\end{figure}

  The orderings of the edges are induced by the counterclockwise orientation of the plane.  Notice that $\Gamma_1$ thickens to a surface of genus zero with four boundary components.  $\G_2$ thickens to a surface of genus 1 with two boundary components.
Of course these surfaces are homotopy equivalent, since they are each homotopy equivalent to the same underlying graph.  But their diffeomorphism types are different, and that is encoded by the cyclic ordering of the vertices.  

These examples make it clear that we need to study the combinatorics of fat graphs more carefully. 
For this purpose, for a fat graph $\G$, let $E(\G)$ be the set of edges, and let $\tilde E(\G)$ be the set of oriented edges.  $\tilde E(\G)$ is a $2$-fold cover of $E(G)$.  It has an involution $E \to \bar E$
which represents changing the orientation.  The cyclic orderings at the vertices determines
a partition of $\tilde E(\G)$ in the following way.  Consider  the  example illustrated in figure 10. .

\begin{figure}[ht]
  \centering
  \includegraphics[height=6cm]{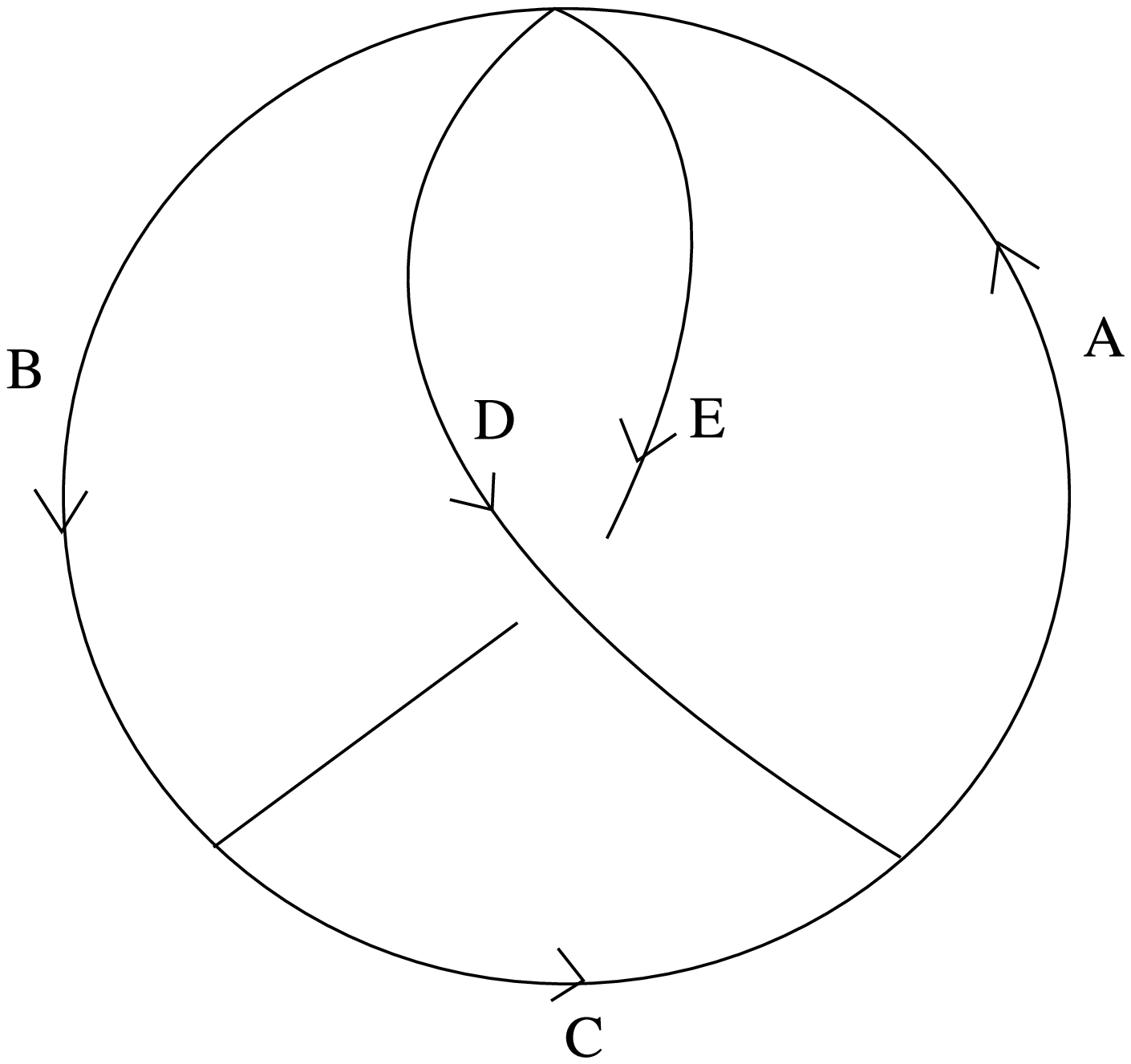}
 \caption{The  fat graph  $\G_2$}
   \label{fig:figten}
\end{figure}
 
 As above, the cyclic orderings at the vertices are determined by the counterclockwise orientation of the plane. 
 To obtain the partition,  notice that an oriented edge has   well defined source and target vertices. Start with an oriented edge, and follow it to its target vertex.  The next edge in the partition is the next oriented edge in the cyclic ordering at that vertex.  Continue in this way until one is back at the original oriented edge.  This will be the first cycle in the partition.  Then continue with this process
 until one exhausts all the oriented edges.  The resulting  cycles in the partition will be called  \sl ``boundary cycles" \rm as they reflect the boundary circles of the thickened surface. In the  case of $\G_2$ illustrated in figure 10, the partition into boundary cycles are given by:
  
  $$
  \text{Boundary cycles of $\G_2$:} \quad (A,B,C) \, (\bar A, \bar D, E, \bar B, D, \bar C, \bar E ).
  $$
  
 So one can compute combinatorially the number of boundary components in the thickened surface of a fat graph.  Furthermore the graph and the surface have the same homotopy type, so one can compute the Euler characteristic of the surface directly from the graph.  Then using the formula
 $\chi (F) = 2 - 2g -n$, where $n$ is the number of boundary components, we can solve for the genus directly in terms of the graph.  The main theorem about spaces of fat graphs is the following (see \cite{penner}, \cite{strebel}).
 
 \begin{theorem}\label{penner}  For $g \geq 2$, the space of metric fat graphs $Fat_{g,n}$ of genus $g$ and $n$ boundary cycles is homotopy equivalent to the moduli space $\cm_{g}^n$ of closed Riemann surfaces of genus $g$ with $n$ marked points.
 \end{theorem} 
 
 Notice that the boundary cycles of a metric fat graph $\G$ nearly determines a parameterization of the boundary of the thickened surface.  For example, the boundary cycle $(A, B, C)$ of the graph $\G_2$ can be represented by a map $S^1 \to \G_2$ where the circle is of circumference equal to the sum of the lengths of sides $A$, $B$, and $C$.  The ambiguity of the parameterization is the choice of where to send the basepoint $1 \in S^1$.  In her thesis \cite{godin}, Godin described the notion of a ``marked" fat graph, and proved the following analogue of theorem \ref{penner}
 
 \begin{theorem}\label{marked}.  Let $Fat^*_{g,n}$ be the space of marked metric fat graphs of genus $g$ and $n$ boundary components.  Then there is a homotopy equivalence
 $$
 Fat^*_{g,n} \simeq \cm_{g,n}
 $$
 where $\cm_{g,n}$ is the moduli space of Riemann surfaces of genus $g$ having $n$ parameterized boundary components.
 \end{theorem}
 
 In \cite{cohengodin} the umkehr map $\rho_{in} : h_*((LM)^p) \to h_{*+\chi (F)\cdot d}(Map (F, M))$
 was constructed as follows.  Let $\G$ be a marked fat graph representing a surface $F$.   Assume $p$ of the boundary cycles of $\G$ have been distinguished as ``incoming", and the remaining $q$ have been distinguished as ``outgoing".  Assume furthermore that  $\G$ satisfies the following technical condition:
 
 \med
 
\begin{definition}\label{chord} A fat graph $\G$ is called a ``Sullivan chord diagram" if it satisfies the following property.  An oriented edge  $E$ is   contained in an incoming boundary cycle of $\G$ if and only if   $\bar E$ is contained in an outgoing boundary cycle. 
 \end{definition}
 
 It is easy to see that every surface $F$ is represented by a marked chord diagram $\G$.  In this case the map
 $$
 \rho_{in} : Map (F, M) \to (LM)^p
 $$
 is homotopic to a map
 $$
 \rho_{in} : Map (\G, M) \to (LM)^p
 $$
 which is obtained by restricting a map from $\G$ to its $p$ incoming boundary cycles, using the parameterizations determined by the markings.   Furthermore it was shown in \cite{cohengodin}
 that this map is a codimension $\chi (F) \cdot d$ embedding, and that a Thom collapse map  could be defined,
 $$
 \tau_F : (LM)^p \to Map (\G, M)^\nu
 $$
 where $\nu$ is the normal bundle of $\rho_{in}$.  This bundle was computed explicitly in \cite{cohengodin}.  This allows for the definition of the umkehr map, as was discussed in the last section.  This in turn, allowed for the definition of the string topology operation $\mu_F$ described above.

 \bg
 We end this section with a discussion of ``open-closed", or perhaps a better term is ``relative"  string topology.   In this setting our background manifold comes equipped with a collection of submanifolds,
 $$
 \scrb = \{ D_i \subset M \}.
 $$
 Such a collection is referred to as a set of ``D-branes", which in string theory supplies
 boundary conditions  for open strings.  In string topology, this is reflected by
 considering the path spaces
 $$
 \cp_M(D_i, D_j) = \{ \gamma : [0,1] \to M, \,: \,  \gamma (0) \in D_i, \quad \gamma (1) \in D_j\}.
 $$
Following \cite{segal}, in a theory with $D$-branes, one associates to a connected, oriented compact one-manifold $S$  whose boundary components are labelled by $D$-branes, a vector space $V_S$.  In the case of string topology,
if $S$ is topologically a  circle, the vector space $V_S = h_*(LM)$.  If $S$ is an interval with boundary points labeled
by $D_i$ and $D_j$, then $V_S = h_*(\cp_M(D_i, D_j))$.  As is usual in field theories, to a disjoint
union of such compact one manifolds, one associates the tensor product of the above vector spaces. 

Now to an appropriate cobordism, one needs to associate an operator between the vector spaces associated to the incoming and outgoing parts of the boundary.  In the presence of $D$-branes these cobordisms are cobordisms of manifolds with boundary.  More precisely, in a theory with $D$-branes, the boundary of a cobordism $F$ is partitioned into three parts:

\begin{enumerate}\label{cobordism}
\item incoming circles and intervals, written $\p_{in}(F)$,
\item outgoing circles and intervals, written $\p_{out}(F)$,
\item the ``free part" of the boundary, written $\p_{f}(F)$, each component of which is labeled by a $D$-brane.  Furthermore $\p_f(F)$ is a cobordism from the boundary of the incoming one manifold to the boundary of  the outgoing one manifold.   This cobordism respects the labeling.
\end{enumerate}  

\begin{figure}[ht]
  \centering
  \includegraphics[height=8cm]{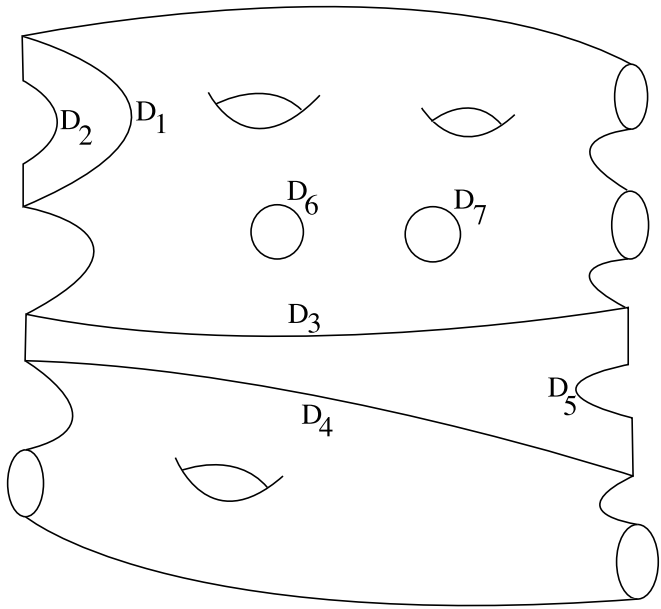}
 \caption{open-closed cobordism }
   \label{fig:figeleven}
\end{figure}

We will call such a cobordism an ``open-closed cobordism" (see figure \ref{fig:figeleven}).  

 In a theory with $D$-branes, associated to such an open-closed cobordism $F$ is an operator,
 $$\phi_F : V_{\p_{in}(F)} \la V_{\p_{out}(F)}.$$  Of course such a theory must respect gluing of open-closed cobordisms.  
 
 Such a theory with $D$-branes  has been put into the  categorical language of PROPs by Ramirez \cite{ramirez} extending notions of Segal and Moore \cite{segal}.  He called such a field theory a $\scrb$-topological quantum field theory.
 
 In the setting of string topology, operators $\phi_F$ were defined by Sullivan \cite{sullivan2} using transversal intersections of chains.  They were defined via Thom-collapse maps by Ramirez in \cite{ramirez} where he proved the following.
 
 \med
 \begin{theorem} Given a set of $D$-branes $\scrb$ in a manifold $M$ and a generalized homology theory $h_*$ that supports orientations of $M$ and all the submanifolds of $\scrb$, then the open-closed string topology operations define a positive boundary $\scrb$- topological quantum field theory.
 \end{theorem}
 
 \section{A Morse theoretic view of string topology}
 The goal of this section is to apply the methods for constructing homology operations using graphs
 described in section 1, to the loop space $LM$, and to thereby recover the string topology operations
 from this Morse theoretic perspective.  This is an exposition of the work to be contained in \cite{cohen2}.  
 
  In order to do this, we need a plentiful supply of Morse functions on $LM$.  Inspired by the work of Salamon and Weber \cite{salamonweber} we take as our Morse functions certain classical energy functions.  These are defined as follows.  Endow our closed $d$-dimensional manifold with a Riemannian metric $g$.  
 
 Consider a potential function on $M$, defined to be a smooth map
 $$
 V : \br/\bz \times M \la \br.
 $$ 
 We can then define the classical energy functional
 \begin{align}\label{energy}
 \cs_V :  LM &\la \br \notag \\
 \gamma &\la \int_0^1\left(\frac{1}{2} |\frac{d\gamma}{dt}|^2 - V(t, \gamma (t)) \right) dt.
 \end{align}
 
 For a generic choice of $V$,  $\cs_V$ is a Morse function \cite{weber}. 
 Its critical points are those $\gamma  \in LM$ satisfying the ODE
 \begin{equation}\label{critical}
 \nabla_t\frac{d\gamma}{dt} = -\nabla V_t(x)
 \end{equation}
 where $\nabla V_t(x)$ is the gradient of the function $V_t (x) = V(t,x)$, and $ \nabla_t\frac{d\gamma}{dt} $ is the Levi-Civita covariant derivative.

  By perturbing $\cs_V$ in a precise way  as  in  \cite{salamonweber}, section 2 , it is possible to assume that $\cs_V$ satisfies the Morse-Smale transversality condition, and one obtains a Morse  chain complex,  $C_*^V(LM)$, for computing the homology of the loop space,
  $H_*(LM)$:
  \begin{equation}\label{morsecomp}
  \to \cdots C_q^V(LM) \xr{\p}C^V_{q-1}(LM)\xr{\p} \cdots
  \end{equation}
  where, as usual, the boundary map is computed by counting gradient trajectories connecting critical points. 
  
  Now as in section 1, we wish to study graph flows, but now the target is the loop
  space, rather than a compact manifold.   Recall that a graph flow is made of gradient trajectories of different Morse functions, that fit together according to the combinatorics
  and the metric of a graph.   In the case of the loop space, a gradient trajectory,
  being a curve in the loop space, may be thought of as a map of a cylinder to $M$.
  In order to fit these cylinders together, we use the combinatorics of a fat graph, as described in section 2.  This is done by the following construction.
  
  \med
  Let $\G$ be a metric marked chord diagram as described in definition \ref{chord}.   Recall this means that the boundary cycles of $\G$ are partitioned into $p$ incoming and $q$ outgoing cycles, and there are parameterizations determined by the markings,
  $$
  \alpha^- : \coprod_p S^1 \la \G, \quad \alpha^+ : \coprod_q S^1 \la \G.
  $$
  
  By taking the circles to have circumference equal to the sum of the lengths of the edges making up the boundary cycle it parameterizes, each component of $\alpha^+$ and $\alpha^-$ is a local isometry.   
  
  Define the surface $\Sigma_\G$ to be the mapping cylinder of these parameterizations,
  \begin{equation}\label{sigmag}
  \Sigma_\G = \left( \coprod_p S^1 \times (-\infty, 0]\right)   \sqcup \left(\coprod_q S^1 \times [0, +\infty)\right)  \bigcup  \G / \sim
  \end{equation}
  where $(t,0) \in S^1 \times (-\infty, 0] \, \sim \alpha^- (t) \in \G$, and $(t,0) \in S^1 \times [0, +\infty)  \, \sim \alpha^+ (t) \in \G$
  
  \med
  
\begin{figure}[ht]
  \centering
  \includegraphics[height=8cm]{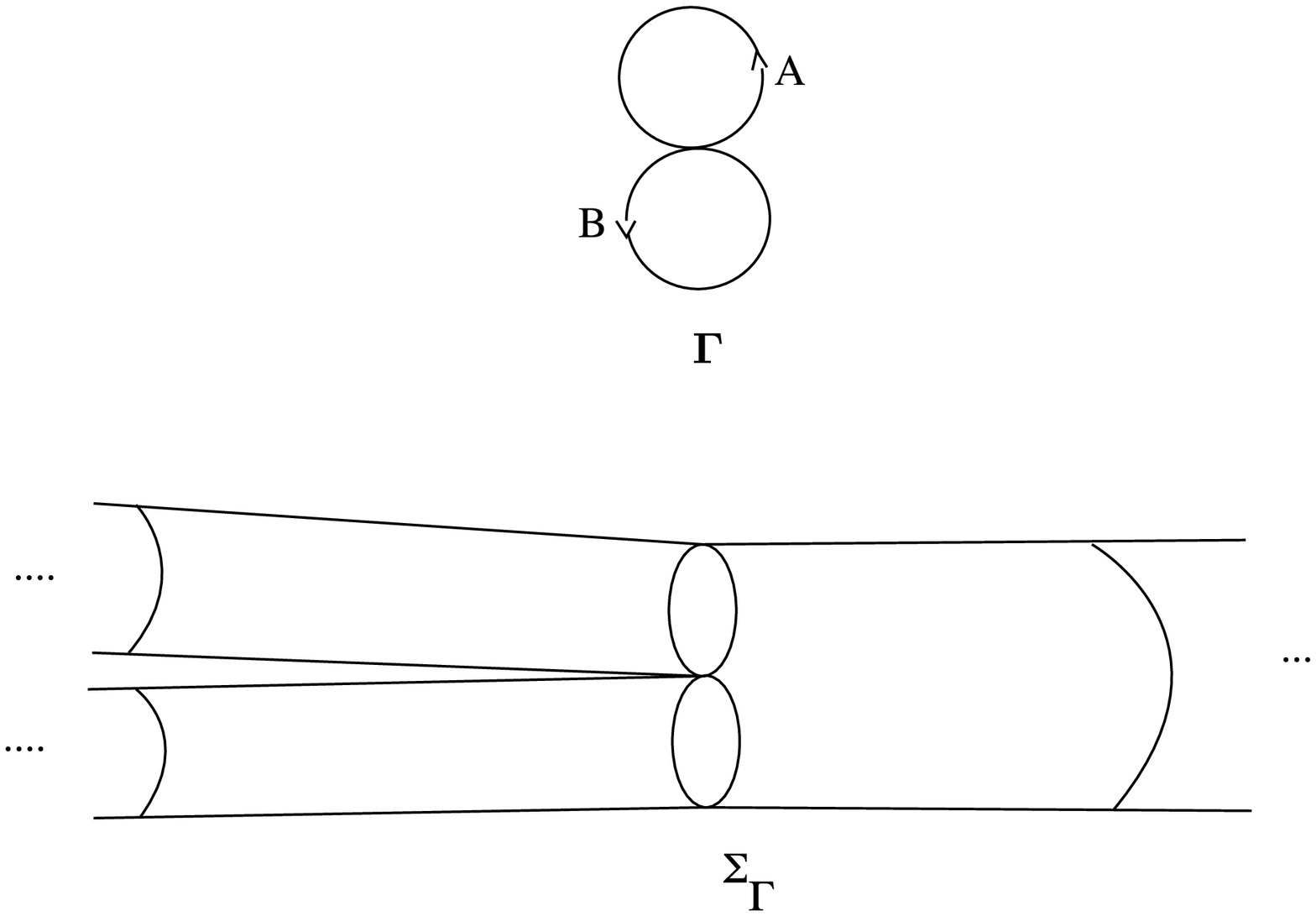}
 \caption{$\Sigma_\G$ }
   \label{fig:figtwelve}
\end{figure}

  Notice that the figure 8  is a  fat graph representing a surface of genus $g=0$ and $3$ boundary components.  This graph has two edges, say $A$ and $B$, and has boundary cycles $(A), (B), (\bar A, \bar B)$.  If we let $(A)$ and $(B)$ be the incoming cycles and $(\bar A, \bar B)$ the outgoing cycle, then the figure 8 graph becomes a chord diagram.
Figure \ref{fig:figtwelve}  is a picture of the surface $\Sigma_\G$, for $\G$ equal to the figure 8.

Notice that a map    $\phi : \Sigma_\G \to M$ is  a collection of $p$ curves  $\phi_i : (-\infty, 0] \to LM$ and $q$-curves,
$\phi_j : [0, +\infty) \to LM$, that have an intersection property at $t=0$ determined by the combinatorics of the fat  graph $\G$. In particular we think of the mapping cylinder $\Sigma_\G$ shown in figure \ref{fig:figtwelve} as the analogue of the ``Y - graph"  (see figure \ref{fig:figfive})
with the vertex ``blown up" via the graph $\G$. 

Now in section 1, a structure on a graph consisted of a metric and a Morse labeling, which was a labeling of the edges by Morse functions.
For example, the Y-graph would have three distinct functions labeling its  three edges.   The analogue of the edges of the Y-graph in our situation are the boundary cylinders of $\Sigma_\G$.  So we need to label these cylinders by functions on the loop space $LM \to \br$.
We choose to restrict our attention to the energy functionals, $\cs_V : LM \to \br$ defined by a potential $V : \br/\bz \times M \to \br$.  This leads to the following defintion.

\begin{definition}  Given a marked chord diagram $\G$ with $p$-incoming and $q$-outgoing boundary cycles, we define
an $LM$-Morse structure $\sigma$  on $\G$ to be a metric  on $\G$ together with a labeling of each boundary cylinder of the surface $\Sigma_\G$ by a distinct potential function $V : \br/\bz \times M \to \br$. (Equivalently, the potential functions label the boundary cycles of the fat graph $\G$.)  
\end{definition}

This leads to the following definition of the moduli space of cylindrical flows.  Compare to Definition \ref{graphflow}.

\begin{definition}  Let $\G$ be a marked chord diagram as above.  Let $\sigma $ be a $LM$-Morse structure on $\G$.  Suppose 
$\phi : \Sigma_\G \to M$ is a continuous map, smooth on the cylinders.   Let $\phi_i : S^1 \times (-\infty, 0] \to M  \, $ be the  restriction of $\phi$ to the $i^{th}$ incoming cylinder, $ i = 1, \cdots, p$,  and $\phi_j : S^1 \times [0, +\infty) \to M  \,$ be the restriction to the $j^{th}$ outgoing cylinder,  $ j  =  1, \cdots, q$. We consider the $\phi_i$'s and $\phi_j$'s as curves in the loop space, $LM$.  Then the moduli space of cylindrical flows  is defined to be
\begin{align}
 \cm^\sigma_\G(LM) = \{\phi : \Sigma_\G \to M \, : \,  \frac{d\phi_i}{dt} + \nabla S_{V_i} (\phi_i(t)) = 0 \quad &\text{and} \quad   \frac{d\phi_j}{dt} + \nabla S_{V_j} (\phi_i(t)) =0 \notag \\
  &\text{for} \quad  i = 1, \cdots , p \quad \text{and} \quad  j = 1, \cdots, q . \} \notag
\end{align}
\end{definition}

\med
  Let $\phi \in \cm^\sigma_\G(LM) $.  For $i = 1, \cdots , p$, let $\phi_{i,-1} : S^1 \to M$ be the restriction of $\phi_i : S^1 \times (-\infty, 0] \to M$ to $S^1 \times \{-1\}$.  Similarly, for $j = 1, \cdots, q$,
  let $\phi_{j,1} : S^1 \to M$ be the restriction of $\phi_j$ to $S^1 \times \{1\}$.  These restrictions define  the following maps (compare with the restriction  and evalation  maps (\ref{eval})  and (\ref{restrict1})).
  
  \begin{equation}\label{restrict2}
 (LM)^q \xleftarrow{\rho_{out}}  \cm^\sigma_\G(LM) \xr{\rho_{in}} (LM)^p .
  \end{equation}
  
  \med
  In \cite{cohen2} it is shown that one can define a Thom collapse map,
  $$
 \tau_\G:  (LM)^p \to (\cm^\sigma_\G(LM)^\nu
 $$
 where $\nu$ is a certain vector bundle of dimension $-  \chi (\G)\cdot d$. This can be thought of as a normal bundle in an appropriate sense.    This allows the definition of an 
 ``umkehr map"
 $$
 (\rho_{in})_! : h_*((LM)^p) \to h_{*+ \chi (\G) \cdot d}(\cm^\sigma_\G(LM))
 $$
 for any homology theory $h_*$ supporting an orientation of $M$.  One can
 then define an operation
\begin{equation}\label{morseop}
 q_\G^{morse} : h_*((LM)^p) \xr{ (\rho_{in})_!}  h_{*+ \chi (\G) \cdot d}(\cm^\sigma_\G(LM)) \xr{(\rho_{out})_*} h_{*+\chi (\G)\cdot d}((LM)^q).
\end{equation}

 The definition of the Thom collapse map $\tau_G$, the induced umkehr map $(\rho_{in})_! $ is a  consequence of the following technical result, proved in \cite{cohen2}.
 
 Consider the map $\psi :  \cm^\sigma_\G(LM) \to Map (\G, M)$ defined by restricting
 a cylindrical flow $\phi : \Sigma_\G \to M$ to the graph $\G \subset \Sigma_\G$.
 
 \begin{theorem}  For a generic choice of $LM$ - Morse structure $\sigma$
   on a marked chord diagram $\G$, the map
   $$
   \psi :  \cm^\sigma_\G(LM) \to Map (\G, M) \simeq Map(\Sigma_\G, M)
      $$
   is a homotopy equivalence.
   \end{theorem}
   
   This result then allows the construction of the Thom collapse map $\tau_\G$ so that the induced umkehr map $(\rho_{in})_!$ is equal to the umkehr map (\ref{restrict1}) via the homotopy equivalence $\psi$.  As a result we have the following, proved in \cite{cohen2}.
   
   \begin{theorem}\label{morsestring} For any marked chord diagram $\G$, the Morse theoretic operation 
   $$q_\G^{morse} : h_*((LM)^p) \la  h_{*+\chi (\G)\cdot d}((LM)^q)$$
   given in (\ref{morseop}) is equal to the string topology operation
   $$
q_\G  : h_*((LM)^p) \la  h_{*+\chi (\G)\cdot d}((LM)^q)
$$
defined  in theorem \ref{stringop}.
\end{theorem}  

This Morse theoretic viewpoint of the string topology operations has another, more
geometric due to Ramirez \cite{ramirez}.  
It is a  direct analogue of the perspective on the graph
operations in \cite{betzcohen}.   

As above, let $\G$ be a marked chord diagram. In Ramirez's setting, an $LM$-Morse structure on $\G$  can involve a labeling of the boundary cycles of $\G$ (or equivalently
the cylinders of $\Sigma_\G$)  by any distinct Morse functions on $LM$ that are bounded below, and  satisfy the Palais-Smale condition  as well as the Morse-Smale transversality condition.   

Let $\sigma$ be an $LM$-Morse structure on $\G$ in this sense.  Let $ (f_1, \cdots, f_{p+q})$  be the Morse functions on $LM$ labeling the $p+q$   cylinders of $\Sigma_\G$.  As above, the first $p$ of these cylinders are incoming, and the remaining $q$ are outgoing.

Let $\vec{a} = (a_1, \cdots , a_{p+q})$ be a sequence of loops such that $a_i  \in LM$
is a critical point of $f_i : LM \to \br$.  Let $W^u(a_i)$ and $W^s(a_i)$ be the unstable and stable manifolds of these critical points. Then define
 
 \med
$\cm^\sigma_\G (LM, \vec{a}) = \{  
  \phi : \Sigma_\G \to M $  that satisfy the following two conditions:   
  \begin{enumerate} 
  \item $\frac{d\phi_i}{dt} + \nabla f_i (\phi_i(t)) = 0 $     for   $  i = 1, \cdots , p+q$  \\
  \item $ \phi_i \in W^u(a_i)$ for $ i = 1, \cdots , p$, and $\phi_j \in W^s(a_j)$ for $j = p+1, \cdots p+q$.\}
  \end{enumerate}
  
  Ramirez then proved that under sufficient transversality conditions described in \cite{ramirez} then $\cm^\sigma_\G (LM, \vec{a})$ is a smooth manifold of dimension
\begin{equation}\label{dimension}
  dim (\cm^\sigma_\G (LM, \vec{a})) = \sum_{i=1}^p Ind (a_i) - \sum_{j=p+1}^{p+q} Ind (a_j)   + \chi (\G)\cdot d.
 \end{equation}
 
 Moreover, an orientation on $M$ induces an orientation on $\cm^\sigma_\G (LM, \vec{a})$.  Furthermore  compactness issues are addressed, and it is shown
 that if $dim (\cm^\sigma_\G (LM, \vec{a})= 0$ then it is compact.  This leads to the following definition.
 For $f_i$ one of the labeling Morse functions, let $C_*^{f_i}(LM)$ be the 
Morse chain complex for computing $H_*(LM)$, and let $ C^*_{f_i}(LM)$ be the corresponding cochain complex.
Consider the chain

\begin{equation}\label{chain}
q_\G^{morse}(LM) =  \sum_{dim (\cm^\sigma_\G (LM, \vec{a}))= 0}  \#\cm^\sigma_\G (LM, \vec{a}) \cdot [\vec{a}]  \quad 
  \in \quad  \bigotimes_{i=1}^p C^*_{f_i}(LM)  \otimes \bigotimes_{j=p+1}^{p+q}C_*^{f_j}(LM)   \end{equation}
  
  We remark that the (co)chain complexes $C^*_{f_i}(LM)$ are generated by critical points, so this large tensor product of chain complexes is generated by vectors of critical points $[\vec{a}]$.     It is shown in \cite{ramirez} that this chain is a cycle and if one uses (arbitrary) field coefficients this defines   a class 
\begin{align}\label{geodef}
 q_\G^{morse}(LM)& \in   (H^*(LM))^{\otimes p} \otimes (H_*(LM))^{\otimes q}  \\
 &= Hom ((H_*(LM))^{\otimes p}, (H_*(LM))^{\otimes q}). \notag
 \end{align}

 Ramirez then proved that 
  these operations are the same as those defined by (\ref{morseop}), and hence by theorem \ref{morsestring}  is equal to the string topology operation.  In the case when $\G$ is the figure 8, then this operation is the same as that defined by Abbondandolo and Schwarz \cite{abschwarz} in the Morse homology of the loop space.

   \section{Cylindrical holomorphic curves in the cotangent bundle}  

This is a somewhat speculative section.  Its goal is to indicate possible relations between string topology operations
and holomorphic curves in the cotangent bundles.  It is motivated by the work of Salamon and Weber \cite{salamonweber}.

As before, we let $M$ be a $d$-dimensional, closed oriented manifold, and $T^*M$ its cotangent bundle.  This is a $2d$-dimensional open manifold with a canonical symplectic form $\omega$ defined as follows.

Let $p : T^*M \to M$ be the projection map.
Let $x\in M$ and $u \in T^*_xM$.  Consider the composition
$$
\alpha (x,u) : T_{(x,u)}(T^*(M)) \xr{Dp} T_x M \xr{u} \br
$$
where $T_{(x,u)}(T^*(M))$ is the tangent space of $T^*(M)$ at $(x,u)$, and $Dp$ is the derivative of $p$. Notice that $\alpha$ is a one form,
$ 
\alpha \in \Omega^1(T^*(M)),
$ 
and we define 
$$
\omega = d\alpha \in \Omega^2(T^*(M)).
$$
It is well known that $\omega$ is a nondegenerate symplectic form on $T^*(M)$.  Now given a Riemannian metric on $M$,
$g : TM \xr{\cong} T^*M$, one gets a corresponding almost complex structure $J_g$ on $T^*(M)$ defined as follows.

\med
The Levi-Civita connection defines a splitting of the tangent bundle of the $T^*(M)$,
$$
T(T^*(M)) \cong p^*(TM) \oplus p^*(T^*(M)).
$$
With respect to this splitting,  $J_g : T(T^*(M)) \to T(T^*(M))$ is defined by the matrix,
$$
J_g = \left( \begin{matrix}
0 & -g^{-1} \\
g & 0
\end{matrix}\right).
$$
The induced metric on $T^*(M)$   is defined by
$$
G_g = \left( \begin{matrix}
g & 0 \\
0 &  g^{-1}
\end{matrix}\right).
$$

 The symplectic action functional is defined on the loop space of the cotangent bundle, $L(T^*(M))$.  Such a loop is given by a pair,
 $(\gamma, \eta)$, where $\gamma : S^1 \to M$, and $\eta (t) \in T^*_{\gamma (t)}M$.  The symplectic action has the formula
 \begin{align}\label{action}
 \ca : L(T^*M) &\to \br  \\
(\gamma, \eta) &\to \int_0^1 \langle \eta (t), \frac{d\gamma}{dt}(t) \rangle dt.
\end{align}

\med
As done by Viterbo \cite{viterbo} and Salamon-Weber \cite{salamonweber},one can do Floer theory on $T^*M$, by perturbing 
the symplectic action functional by a Hamiltonian induced by a potential function $V : \br/\bz \times M \to \br$ in  the following way.

Given such a potential $V$, define $H_V : \br/\bz \times T^*(M) \to \br$ by the formula
\begin{equation}\label{hamilton}
H_V(t, (x,u)) = \frac{1}{2}|u|^2 + V(t, x).
\end{equation}
Then one has a perturbed symplectic action
\begin{align}\label{AV}
\ca_V : L(T^*M) &\to \br \\
(\gamma, \eta) &\to \ca (\gamma, \eta) - \int_0^1H(t, (\gamma (t), \eta (t)))dt.
\end{align}

As observed in \cite{salamonweber}, via the Legendre transform one sees that the critical points of $\ca_V$ are loops
$(\gamma, \eta)$, where $\gamma \in LM$ is a critical point of the energy functional $\cs_V : LM \to \br$, and $\eta$ is determined
by the derivative $\frac{d\gamma}{dt}$ via the metric, $\eta (v) = \langle v,\frac{d\gamma}{dt}\rangle$.
Thus the critical points of $\ca_V$ and those of $\cs_V$ are in bijective correspondence.  These generate the Floer complex,
$CF^{ V}_*(T^*M)$ and the Morse complex, $C^V_*(LM)$ respectively.  The following result is stated in a form proved by by Salamon and Weber in \cite{salamonweber}, but the conclusion of the theorem was first proved by Viterbo \cite{viterbo}.

\begin{theorem}\label{floermorse} The Floer chain complex  $CF^{ V}_*(T^*M)$ and the Morse complex $C^V_*(LM)$
are chain homotopy equivalent.  There is a resulting isomorphism of the Floer homology of the cotangent bundle with the 
homology of the loop space,
$$
HF^{V}_*(T^*M) \cong H_*(LM).
$$
\end{theorem}

This result was also proved using somewhat different methods by Abbondandolo and Schwarz \cite{abschwarz}.

The Salamon-Weber argument involved scaling the metric on $M$,  $g \to \frac{1}{\eps}g$, which  scales the almost complex
structure $J \to J_\eps$,  and the metric on $T^*M$, $G \to G_\eps =   \left( \begin{matrix}
\frac{1}{\eps}g & 0 \\
0 &  \eps g^{-1}
\end{matrix}\right).$ Notice that in this metric, the ``vertical" distance in the cotangent space is scaled by $\eps$.  

Now the boundary operator in the Floer complex $CF^{V}_*(T^*M)$ is defined by counting gradient flow lines of $\ca_V$, which
are curves $(u,v): \br \to  \ L(T^*M)$, or equivalently,
$$
(u,v) : \br \times S^1 \to T^*M
$$ that satisfy the perturbed Cauchy Riemann equations,
\begin{equation}\label{cauchyriemann}
\p_su - \nabla_t v - \nabla V_t (u) = 0   \quad \text{and} \quad \nabla_s v +\p_t u - v = 0.
\end{equation}

We refer to these maps as   holomorphic cylinders in $T^*M$ with respect to the almost complex structure $J_\eps$ and the Hamiltonian
$H_V$.    Salamon and Weber proved that there is an $\eps_0 >0$ so that for $\eps < \eps_0$,  the set of these  holomorphic cylinders     defined with respect to the metric $G_\eps$, that  connect  critical points $(a_1, b_1)$ and  $(a_2, b_2)$ of relative Conley-Zehnder index one, is in bijective correspondence with the set of gradient 
trajectories of the energy functional $\cs_V : LM \to \br$ defined with respect to the metric $\frac{1}{\eps}g$ that connect $a_1$ to $a_2$.
Theorem \ref{floermorse} is then a consequence. 

The Salamon-Weber construction inspires  the following idea.  Let $\G$ be a marked chord diagram as in the previous section.
Let $\Sigma_G$ be the cylindrical surface, and let $\sigma$ be an $LM$ - Morse structure on $\G$.  Recall that this means
that the graph $\G$ has a metric, and hence the cylinders $S^1 \times [0, +\infty)$ and $S^1 \times (-\infty, 0]$ making up $\Sigma_\G$ have well defined widths, and hence complex structures. Recall that part of the data of a $LM$-Morse structure $\sigma$ is a labeling of these boundary cylinders by distinct potentials, $V_i : \br/\bz \times M \to \br$.  Given a Riemannian metric $g$ on $M$ as above, and an $\eps > 0$, we define the moduli space of  ``cylindrical holomorphic curves" in the cotangent bundle $T^*(M)$ as follows.

\begin{definition}\label{holo}
We define $\cm^{hol}_{(\G, \sigma, \eps)}(T^*M)$ to be the space of continuous maps
$  \phi : \Sigma_\G \to T^*(M)$  such that the restrictions to the cylinders, 
$$
\phi_i : (-\infty, 0] \times S^1_{c_i} \to T^*M \quad \text{and} \quad \phi_j : [0, +\infty) \times S^1_{c_j} \to T^*M
$$
are holomorphic with respect to the almost complex structure $J_{\eps}$ and the Hamiltonians $H_{V_i}$ and $H_{V_j}$ respectively.
Here the circles $S^1_{c}$ are round with circumference $c_j$ determined by the metric given by the structure $\sigma$.
\end{definition}  

\med
Like in the last section we have restriction  maps (compare (\ref{restrict2}))
 \begin{equation}\label{restrict3}
 (L(T^*M))^q \xleftarrow{\rho_{out}} \cm^{hol}_{(\G, \sigma, \eps)}(T^*M) \xr{\rho_{in}} (L(T^*M))^p .
  \end{equation}
 $\rho_{in}$ is  defined
by  sending a cylindrical flow $\phi$ to  $\prod_{i=1}^p \phi_{i,-1} : \{-1\} \times S^1  \to T^*M$ and $\rho_{out}$ sends $\phi$ to $\prod_{j=p+1}^{p+q} \phi_{j,1}:  \{1\} \times  S^1\to T^*M.$

We conjecture the following analogue of the existence of the string topology operations, and their field theoretic properties.  

\begin{conjecture}\label{basic} For every marked chord diagram $\G$, there is an umkehr map
$$
(\rho_{in})_! : (HF_*(T^*M))^{\otimes p} \to H_{*+\chi (\G)\cdot d}( \cm^{hol}_{(\G, \sigma, \eps)}(T^*M))
$$
and a homomorphism
$$
(\rho_{out})_* :  H_{* }( \cm^{hol}_{(\G, \sigma, \eps)}(T^*M)) \to (HF_*(T^*M))^{\otimes q}
$$
so that the operations
$$
\theta_\G  = (\rho_{out})_* \circ (\rho_{in})_! :  (HF_*(T^*M))^{\otimes p}  \to (HF_*(T^*M))^{\otimes q}
$$
satisfy the following properties:
\begin{enumerate}
\item The maps $\theta$ fit together to define a positive boundary, topological field theory.
\item With respect to the Salamon-Weber isomorphism $HF_*(T^*M) \cong H_*(LM)$ (theorem \ref{floermorse})
the Floer theory operations $\theta_\G$ equal the string topology operations $q_\G$ studied in the last two sections.
\end{enumerate}
\end{conjecture}

\bfl
\bf Remark. \rm
   The existence of a field theory structure on the Floer homology of a closed symplectic manifold was established by Lalonde \cite{lalonde}.  The above conjecture should be directly related to Lalonde's constructions.
   
   \efl
  \med
  A possible way to approach this conjecture is to prove the following, which can be viewed as a generalization
  of the Salamon-Weber result relating the gradient trajectories of the Floer functional $\ca_V$ on $L(T^*M)$ (i.e $J$-holomorphic cylinders in $T^*M$)  with
  the gradient trajectories of the energy functional $S_V$ on $LM$.

 \begin{conjecture}
 There is a natural map induced by the Legendre transform,
 $$
 \ell : \cm^\sigma_\G(LM))   \to   \cm^{hol}_{(\G, \sigma, \eps)}(T^*M))   
 $$
 which is a homotopy equivalence of $\eps$ sufficiently small.
 \end{conjecture}
 
In view of theorem \ref{morsestring},  this conjecture would imply conjecture \ref{basic} .
 
 \med
 Now one might also take the more geometric approach to the construction of these Floer theoretic operations, analogous to Ramirez's  geometrically defined Morse  theoretic constructions of string topology operations.  This would involve the study of the space of cylindrical holomorphic curves in $T^*M$, with boundary conditions   in stable
 and unstable manifolds of critical points, $ \cm^{hol}_{(\G, \sigma, \eps)}(T^*M, \vec{a})$. Smoothness and compactness properties need to be established for these moduli spaces.  In particular, in a generic situation  their dimensions should
 be given by the formula
 $$
 dim \,( \cm^{hol}_{(\G, \sigma, \eps)}(T^*M, \vec{a})) = \sum_{i=1}^p Ind (a_i)  - \sum_{j=p+1}^{p+q}Ind (a_j)  + \chi (\G)\cdot d
 $$
 where $Ind (a_i)$ denotes the Conley-Zehnder index. 
 
 \med
 We remark that in the case of the figure 8,  this analysis has all been worked out by Abbondandolo and Schwarz \cite{abschwarz2}.  In this case $\Sigma_\G$ is a Riemann surface structure on the pair of pants.  They proved the existence of a ``pair of pants" algebra structure on $HF_*^V(LM)$ and with respect to their isomorphism,
 $HF_*^V(LM) \cong H_*(LM)$ it is isomorphic to the pair of pants product on the Morse homology of $LM$.  In view
 of the comment following definition \ref{geodef} we have the following consequence.
 
 \med
 \begin{theorem} With respect to the isomorphism $HF^V_*(T^*M) \cong H_*(LM)$, the pair of pants product in the Floer homology of the cotangent bundle corresponds to the Chas-Sullivan string topology product.
 \end{theorem}
 
 \med
 Another aspect of the relationship between the symplectic structure of the cotangent bundle and the string topology
 of the manifold,  has to do with the relationship between the moduli space of $J$-holomorphic curves with cylindrical
 boundaries, $\cm_{g,n}(T^*M)$, and moduli space of cylindrical holomorphic curves, $ \cm^{hol}_{(\G, \sigma, \eps)}(T^*M)$, where we now let $\G$ and $\sigma$ vary over the appropriate space of metric graphs.  These moduli spaces should be related as a parameterized version of the relationship between the moduli space of Riemann surfaces and the space of metric fat graphs (theorem \ref{penner}). 
 
 Once established, this relationship would give a direct relationship between Gromov-Witten invariants of the cotangent bundle, and the string topology of  the underlying manifold. In this setting the Gromov-Witten invariants would be  defined using moduli spaces of curves with cylindrical ends rather than marked points, so that the invariants would be defined in terms of the homology of the loop space (or, equivalently, the Floer homology of the cotangent bundle), rather than the homology of the manifold. 
 
 One might also speculate about the relative invariants.  As described in section 2, there is an ``open-closed" version of string topology, defined in the presence of $D$-branes, which are submanifolds $D_i \subset M$.  Recall that on the cotangent level, the conormal bundles,
 $$
 Conorm (D_i) \subset T^*M
 $$
are Lagrangian submanifolds.  It is interesting to speculate about the open-closed string topology invariants, defined on the homology of the space of paths, $H_*(\cp_M(D_i, D_j))$ and how they may be related to the relative Gromov-Witten invariants of the conormal bundles in the cotangent space, or the similar (  more general) invariants of the cotangent bundle  coming
from the symplectic field theory of Eliashberg,  Givental, and Hofer \cite{egh}.  We believe that the relationship between the symplectic topology of the cotangent bundle and the string topology of the underlying manifold is  very rich.

\end{document}